%% file: main.tex
\documentclass{amsart}
\usepackage{bbm}
\usepackage{graphicx}
\usepackage{verbatim}
\usepackage{stackrel}
\usepackage{color}
\graphicspath{{Fig/}}
\usepackage{amsmath}
\usepackage{amssymb}
\usepackage{array}
\usepackage{mathtools}
\usepackage{xcolor}
\usepackage{mathdots}
\usepackage{pgf}
\usepackage{bm}
\usepackage{hyperref}
\usepackage{CJKutf8}

\usepackage{enumerate}

\usepackage{enumitem}

\usepackage{tikz}
\usetikzlibrary{shapes.arrows, fadings}
\usetikzlibrary{arrows}
\usetikzlibrary{matrix}
\usetikzlibrary{shapes,intersections}

\hypersetup{
    colorlinks,
    linkcolor={blue},
    citecolor={red},
}

\newtheoremstyle%
 {bluethm}%
 {}{}%
 {\color{blue}\itshape}
 {}%
 {\color{blue}\bfseries}%
 {\color{blue}.}%
 { }{}

\newtheoremstyle%
 {greenthm}%
 {}{}%
 {\color{green!50!black!100!}\itshape}
 {}%
 {\color{green!50!black!100!}\bfseries}%
 {\color{green!50!black!100!}.}%
 { }{}

 \newtheoremstyle%
 {redthm}%
 {}{}%
 {\color{red}\itshape}
 {}%
 {\color{red}\bfseries}%
 {\color{red}.}%
 { }{}
\newtheorem{theorem}{Theorem}
\newtheorem{prop}{Proposition}
\newtheorem{lemma}{Lemma}

\newtheorem{corollary}{Corollary}

\newtheorem{example}{Example}

\theoremstyle{definition}

\newtheorem{defi}{Definition}

\newtheorem{remark}{Remark}

\newtheorem{assumption}{Assumption}[part]

\def\N{{\mathbb N}}
\def\R{{\mathbb R}}
\def\bS{{\mathbb S}}
\def\S{{\mathcal S}}

\def\P{{\mathbb P}}
\def\E{{\mathbb E}}

\def\cc{c}
\def\mmu{\mu}
\def\B{{B_\delta(\beta)}}
\def\tX{{\tilde{X}}}
\def\<{{\langle}}
\def\>{{\rangle}}
\def\sumn{{\sum_{i=1}^{N(Z_n n)}}}
\def\mdp{{\frac{1}{a_n^2}\log\P}}
\def\mdpsup{{\limsup_{n\rightarrow \infty}\frac{1}{a_n^2}\log\P}}
\def\mdpinf{{\liminf_{n\rightarrow \infty}\frac{1}{a_n^2}\log\P}}
\def\limn{{\lim_{n\rightarrow\infty}}}
\newcommand{\diff}{\mathop{}\mathopen{}\mathrm{d}}

\newcommand\ind[1]{\mathbbm{1}_{\left\{#1\right\}}}

\def\cal{\mathcal}
\def\eps{\varepsilon}

\makeatletter
\newcommand{\proofstep}[1]{%
  \par
  \addvspace{\medskipamount}
  \textit{#1\@addpunct{.}}\enspace\ignorespaces
}
\makeatother

\title[MDPs in the multi-type random graphs]{On the moderate deviation principles in the  sparse multi-type Erd\H{o}s R\'enyi random graph}

\author{Rui Yu}
\email{yurui42@mail.ustc.edu.cn}

\author{Wen Sun}
\email{wensun.ustc@gmail.com}
\address        {School of Mathematical Sciences, University of Science and Technology of China, Jinzhai 96, 230026 Hefei}

\date{\today}
\keywords{ Multi-type Erd\H{o}s-R\'enyi random graph;  Moderate deviation principle;  Compound Poisson process; Conditional limit theorem; Component size.}

\begin{document}

\maketitle

\input Abstract

\bigskip

\hrule

\vspace{-2mm}

\tableofcontents

\vspace{-5mm}

\hrule

\bigskip

\input intro
\input cpp

\input connect

\input mdpcpp

\input super


\appendix
\input app

\section*{Acknowledgments}
 This work is supported by the National Key R\&D Program of China (No. 2022YFA1006500) and by the National Natural Science Foundation of China (No. 12401171).

\bibliographystyle{amsplain}
\bibliography{ref}
\bigskip

\end{document}

%% file: Abstract.tex
 This paper investigate the sparse multi-type Erd\H{o}s R\'enyi random graphs studied in S\"{o}derberg~\cite{soderberg2002general} and also Bollob\'as et al.~\cite{bollobas2007phase}. Although the corresponding central limit results are currently unknown, we establish moderate deviation principles for the size of the largest connected component, the number of specific types of connected components, and the total number of connected components. The associated rate functions are provided explicitly. As a byproduct of this work, we present the law of large numbers for the total number of connected components. Our proof methodology relies on representing the multi-type random graph using a conditional multi-dimensional compound Poisson process. We also discuss the properties of related multi-type branching processes and the properties of the matrices in the rate functions.

%% file: intro.tex
\section{Introduction and Main results}
The classical sparse Erd\H{o}s-R\'enyi random graph $G(n,t/n)$, which is obtained by adding edges independently with probability $t/n$, for some constant $t>0$, to the vertex set $[n]:=\{1,2,\dots,n\}$, has been a central topic in the study of random graphs since 1950s. It is well-known from the original paper~\cite{MR0125031} by Erd\H{o}s and R\'enyi that a phase transition occurs when $t>1$.  In this super-critical regime, a giant connected component appears and has a size approximately equal to $\left(1-T/t\right)n$, where  $T< 1$ and satisfies a duality relation $Te^{-T}=te^{-t}$. 
The central limit theorem for the size of this giant connected component is proven by  Pittel~\cite{MR1099795} through a study of the number of trees in the graph and  the associated ordinary differential equations, by Martin-L\"of~\cite{MR1659544} through a study of the SIR model and the asymptotic stochastic differential equations and by Barraez et al.~\cite{MR1786919} through an analysis of a depth-first search algorithm. The moderate deviation principle for this giant component is proven by  Puhalskii~\cite{MR2118868} by studying a queuing system, by Sun~\cite{sun2023conditional} by using a compound Poisson process representation. The large deviation principle for this giant component is proved by O'Connell~\cite{MR1616567}.
In addition to the giant component, the empirical measure of all the components and the total number of the components are also well understood. For example, see Pittel~\cite{MR1099795} for the law of large numbers and CLTs, Puhalskii~\cite{MR2118868} for CLTs and LDPs, Sun~\cite{sun2023conditional} for MDPs and Andreis et al.~\cite{MR4323309} for LDPs. See also books~\cite{MR1140703,MR0809996,MR3617364} for other related studies.

In this paper, we study a sparse multi-type Erd\H{o}s and R\'enyi random graph proposed by S\"{o}derberg~\cite{soderberg2002general, soderberg2003properties, soderberg2003random}.
For $n\in\N$, let $G_n=\mathcal{G}(n,\mu^n,\kappa)$ be a random graph defined on the vertices set $[n]=\{1,, \cdots, n\}$, where each vertex is assigned a type from a non-empty and finite set $\S=\{1,\cdots,d\}$, known as the type set. 
Let $\mu^n=(\mu^n_i,i\in\cal{S})$ be a probability measure on $\S$ and $\mu^n_in\in\N$ for all $i\in\S$.
We suppose that the number of vertices of type $r$ is exactly $\mu^n_r n$.
Furthermore, we assume that as $n\rightarrow\infty$, $\mu^n$ converges weakly to some probability measure $\mu=(\mu_r)_{r\in\S}$ with $\mu_r>0$, $\forall r\in\cal{S}$.
The graph is undirected, excluding self-loops and multiple edges. All possible edges are sampled independently, with the probability of an edge between two vertices of types $r$ and $s$ given by $1\wedge \frac{\kappa(r,s)}{n}$, where the kernel $\kappa:\S\times\S\rightarrow(0,\infty)$  is a positive symmetric bounded function. We also use $\kappa$ to represent the matrix $(\kappa(r,s))_{r,s\in\S}$ throughout the paper. For any vertex $v\in G_n$, let $C(v)\in \N^d$ be the type configuration of the connected component  containing $v$, \emph{i.e.} $C(v)_s$ is the number of type $s$ vertices for all $s\in\cal{S}$ in this connected component containing $v$.

\subsection{The phase transition}
It has been shown that when the largest eigenvalue, denoted by $\Sigma(\kappa,\mu)$, of the matrix
\[
\kappa D_\mu:=(\kappa(r,s)\mu_s,r,s\in\cal{S})
\]
is larger than $1$, then a unique giant component appears in the random graph  $G_n$ as $n\to \infty$. Moreover, the type configuration of this giant component, denoted by $\cal{C}^n_{\rm max}$,  is asymptotically equal to $(\mu-c)n$, where $c=(c_i,i\in\S)$  solves the following characteristic equation with respect to  $\mu$,
    \begin{equation}\label{cha}
        c_ie^{-(\kappa c)_i}=\mu_i e^{-(\kappa\mu)_i}, \forall i\in\S,
    \end{equation}
and satisfies $0<c_i<\mu_i$ for all $i$ and $\Sigma(\kappa,c)<1$.
See Theorem~6.2 in Bollob\'as et al.~\cite{bollobas2007phase} and Theorem~1 in  S\"{o}derberg~\cite{soderberg2002general} for more details.

   \subsection{The empirical measure on the small components} The large deviations of the empirical measure of the connected components have been studied by Andreis et al.~\cite{LDPinhomo}. As a byproduct, we know that, for any $k\in\N^{d}$, the number of connected components with configuration $k=(k_s,s\in\cal{S})$, denoted by $t_n(k)$, is approximately $h(k)n$ for $n$ large, where
\begin{equation}\label{hk}
    h(k):=\tau(k)\prod_{s\in\S}\frac{(e^{-(\kappa\mu)_s}\mu_s)^{k_s}}{k_s!},
\end{equation}
   \begin{equation}\label{tau(k)}
    \tau(k):=\sum_{T\in\cal{T}(k)}\prod_{i,j\in E(T)}\kappa(x_i,x_j),
    \end{equation}
and $\cal{T}(k)$ is the set of spanning trees on $[|k|]$, $E(T)$ is the set of edges of tree $T$ and $(x_1,\cdots,x_{|k|})\in\S^{|k|}$ is a vector such that $k_s=\sum_{i=1}^{|k|}\ind{x_i=s}$.

\subsection{The multi-type Galton-Watson tree}\label{branching} The component structure in the random graph is often studied via a branching process that mimics the graph exploration algorithm. In our case, this is a multi-type Poisson Galton-Watson branching process on the type space $\cal{S}$.A particle of type $s$ yields subsequent generations where the number of type-$r$ particles independently follows a Poisson distribution with intensity $\kappa(s,r)\mu_r$. Bollob\'as et al.~\cite{bollobas2007phase} characterized various properties of this process. For example, for a multi-type branching process starting from a type-$r$ particle,  the probability that the total population $Y=(Y_s,s\in\cal{S})$ has exactly $k_s$ type-$s$ particles, $\forall s\in\cal{S}$, is given by
\[
\P_r(Y=k)=\frac{k_r}{c_r}h(k),
\]
where $(c_r,r\in\cal{S})$ is the solution to the characteristic equation~\eqref{cha} and $h(k)$ is given in~\eqref{hk}. Further details on multi-type branching processes can be found in Harris~\cite{harris1963theory} and Kallenberg~\cite{kallenberg1997foundations}.

Similarly to the classical `one-type' Erd\H{o}s and R\'enyi random graph, the distribution of the type vector $Y$ is closely related to the type-configuration of a  connected component in $G_n$. Since all connected components are identically distributed, intuitively, the  total number of the connected components,
\[
C_n:=\sum_{k\in\N^{d}}t_n(k),
\]
should be close to the ratio between the total number of vertices in the ``non-giant'' connected components and the total population of the branching process described above. See the author's prior work~\cite{sun2023conditional} for the `one-type' case. We  prove this intuition in the multi-type graph rigorously later in Corollary~\ref{corcn}, showing that the total number $C_n$ of connected components in $G_n$ is approximately
\[\left(|c|-\frac{1}{2}c^T\kappa c\right)n,\]
for large $n$. We show later in Section~\ref{sec:jumps} that this limit is equal to the expected value $\sum_{r\in\cal{S}}c_r\E_r[1/|Y|]$, where $(c_r,r
\in\S)$ can be seen as the fraction of types in the ``non-giant'' particles.
\subsection{Main results}
Our main purpose is to construct a conditional multi-dimensional compound Poisson process that represents the type configurations of all connected components in the sparse multi-type Erd\H{o}s and R\'enyi random graph $G_n$. We  use this representation to study the fluctuations in the model. More specifically, we will establish the moderate deviation principles for the type configuration  of the giant component $\cal{C}_{\rm max}^n$, the number of a certain type of connected components $t_n(k)$ and the total number of connected components $C_n$. We note that, to the best of our knowledge, the central limit theorems of these three quantities are still unknown. While our method is not sufficient to obtain these central limits, we conjecture that the related fluctuations are Gaussian, with variances specified by our rate functions.

Before our main results, we introduce an assumption on $(\kappa,\mu)$ that we will use throughout the paper. It is a sufficient condition for the convergence of series in Proposition~\ref{moment}. We believe this condition may be stronger than necessary, but it simplifies the proofs by avoiding excessive technicalities. We mention that in our previous paper~\cite{sun2023conditional} on the classical Erd\H{o}s-R\'enyi graph, we overcame this assumption using analytical methods.

\begin{assumption}\label{assump1}
Let $[\kappa]:=\sup\{\<\nu,\kappa\nu\>:\nu\in\R^d, |\nu|=1, \nu_s>0, s\in\S\}$. Suppose $\kappa$, $\mu$ satisfy 
$\Sigma(\kappa,\mu)-\log\Sigma(\kappa,\mu)-\frac12[\kappa]>1.$
\end{assumption}

 We now state the moderate deviations in the supercritical setting, where the giant component appears. In the following, for any vector $u\in\R^d$,  $D_u$ denotes the $d\times d$ diagonal matrix $\textrm{diag}(u)$.
\begin{theorem}\label{thm1}
  Under Assumption \ref{assump1}, when  the random graph $G_n$ is super-critical $(\Sigma(\kappa,\mu)>1)$, for any sequence $(a_n)$ such that $n^{1/4}\ll a_n\ll \sqrt{n}$ and $|\mu^n-\mu|=O(\frac{1}n)$,  the sequence of random vectors
\[
\frac{1}{a_n\sqrt{n}}\left(\cal{C}_{\rm max}^n-(\mu-\cc)n\right)
\]
satisfies the moderate deviation principle in $\R^d$ with speed $a_n^2$ and the rate function
\[
\cal{I}(x)=\frac{1}{2}\left\langle\left(I-D_\cc\kappa\right)x,D_v\left(I-D_\cc\kappa\right)x\right\rangle,
\]
where $D_v$ is the diagonal matrix with the vector
\[
v=\left(\frac{\mu_i}{\cc_i(\mu_i-\cc_i)},i\in\cal{S}\right).
\]
\end{theorem}

\begin{theorem}\label{thm2}
Under Assumption \ref{assump1}, when the random graph $G_n$ is super-critical $(\Sigma(\kappa,\mu)>1)$,  for any fixed $k\in\N^d$, any sequence $(a_n)$ such that $n^{1/4}\ll a_n\ll \sqrt{n}$ and $|\mu^n-\mu|=O(\frac{1}n)$, 
the sequence of random variables
\[
\frac{1}
{a_n\sqrt{n}}\left(t_n(k)-h(k)n
\right)
\]
satisfies the moderate deviation principle in $\R$ with speed $a_n^2$ and the rate function
\[
        J_k(x)=\frac{x^2}{2}\left(\frac{1}{h(k)}+k^TA(A+B_k)^{-1}B_kk\right),
        \]
        where 
$A=(I-\kappa D_\cc )D^{-1}_{\mu-\cc}(I-D_\mu \kappa)$, $\Phi=(D_c^{-1}-\kappa)^{-1}$ and $B_k^{-1}=\Phi-kk^Th(k)$.
\end{theorem}

\begin{theorem}\label{thm3}
Under Assumption \ref{assump1}, when  the random graph $G_n$ is super-critical $(\Sigma(\kappa,\mu)>1)$, for any sequence $(a_n)$ such that $n^{1/4}\ll a_n\ll \sqrt{n}$ and $|\mu^n-\mu|=O(\frac{1}n)$,   the sequence of random variables
\[
\frac{1}
     {a_n\sqrt{n}}
\left(
     C_n-\left(|\cc|-\frac{1}{2}\langle\cc,\kappa\cc\rangle\right)n
\right)
\]
satisfies the moderate deviation principle in $\R$ with speed $a_n^2$ and the rate function
\[
i(x)=\frac{x^2}{2}
    \frac{1}{|\cc|-\frac{1}{2}\langle\cc,\kappa\cc\rangle}+\frac12\left(\frac{x}{|\cc|-\frac 12\langle \cc, \kappa\cc\rangle}\right)^2
    \cc^TA(A+B)^{-1}B\cc,
    \]
where $A$ and $\Phi$ are given in Theorem~\ref{thm2} and $B^{-1}=\Phi-(|\cc|-\frac 12\langle \cc, \kappa\cc\rangle)^{-1}\cc \cc^T$.
\end{theorem}
The moderate deviation principles for the marginal empirical measure $t_{n}(k)$ and the total number of connected components $C_n$ in the subcritical graph $G_{n}$ are also established. Since their corresponding proofs are analogous to those for the super-critical case, we present only the main results here.
\begin{theorem}\label{thm4}
  Under Assumption \ref{assump1}, when the random graph $G_n$ is sub-critical $(\Sigma(\kappa,\mu)<1)$, for any fixed $k\in\N^d$, any sequence $(a_n)$ such that $\sqrt{\log n}\ll a_n\ll \sqrt{n}$ and $|\mu^n-\mu|=O(\frac{1}n)$,   the sequence of random variables
  \[
\frac{1}
{a_n\sqrt{n}}\left(t_n(k)-h(k)n
\right)
\]
satisfies the moderate deviation principle in $\R$ with speed $a_n^2$ and the rate function
\[
\tilde{J}_k(x)=\frac{x^2}{2}\left(\frac{1}{h(k)}+k^T\left(\Phi-kk^Th(k)\right)^{-1}k\right).
\]
\end{theorem}

\begin{theorem}\label{thm5}
  Under Assumption \ref{assump1}, when the random graph $G_n$ is sub-critical $(\Sigma(\kappa,\mu)<1)$, for any sequence $(a_n)$ such that $\sqrt{\log n}\ll a_n\ll \sqrt{n}$ and $|\mu^n-\mu|=O(\frac{1}n)$,   the sequence of random variables
  \[
\frac{1}
     {a_n\sqrt{n}}
\left(
     C_n-\left(1-\frac{1}{2}\langle\mu,\kappa\mu\rangle\right)n
\right)
\]
satisfies the moderate deviation principle in $\R$ with speed $a_n^2$ and the rate function
\[
\tilde{i}(x)=\frac{x^2}{2}\frac{1}{1-\frac12\<\mmu,\kappa\mmu\>}+\frac{1}{2}\left(\frac{x}{1-\frac12\<\mmu,\kappa\mmu\>}\right)^2\mu^T\left((\Phi-(1-\frac{1}{2}\<\mu,\kappa\mu\>)^{-1}\mu\mu^T\right).
\]
\end{theorem}

As an application of our Theorem~\ref{thm3} and Theorem~\ref{thm5}, combined with the Borel-Cantelli lemma, we can immediately derive the law of large numbers for the total number of connected components $C_n$ as follows.

\begin{corollary}[LLN for $C_n$]\label{corcn}
Under Assumption~\ref{assump1}, when the random graph is non-critical ($\Sigma(\kappa,\mu)\neq 1$),  we have, for the almost surely convergence, 
\[
\lim_{n\to\infty}\frac{C_n}{n}= |\cc|-\frac{1}{2}\langle\cc,\kappa\cc\rangle.
\]
\end{corollary}
\subsection*{Outline of the paper}
The rest of the paper is organized as follows. In Section~\ref{sec:cpp}, we define the multi-dimensional compound Poisson process that we use to study the random graph $G_n$, and we also restate our main results in terms of this process. In Section~\ref{sec:connect}, we analyze the asymptotic behavior of the process as $n$ tends to infinity. Section~\ref{sec:mdpcpp} presents the moderate deviation principles for this compound Poisson process, along with the properties of the matrices involved in the rate functions. Finally, in Section~\ref{sec:super}, we prove the main results concerning the super-critical regime of the graph.

\subsection*{Notations}
We denote by $ \N = \{0,1,2,\dots\} $ and by $ \N^+ = \{1,2,3,\dots\} $. For $ n \in \N $, we write $ f(n) = O(g(n)) $ if there exist a positive constant $ M > 0 $ and an integer $ n_0 \in \N $ such that
$$ |f(n)| \leq M g(n), \quad \forall n \geq n_0. $$
For $ n \in \N $ and $ \xi \in \R $, we say that $ f(n,\xi) = O_\xi(g(n,\xi)) $ if, for every compact set $ K \subset \R $, there exist a positive constant $ M > 0 $ and an integer $ n_0 \in \N $ such that
$$ |f(n,\xi)| \leq M g(n,\xi), \quad \forall n \geq n_0, \quad \forall \xi \in K. $$
For any $ x \in \R^d $, we denote by $ D_x $ the diagonal matrix $ \operatorname{Diag}(x) $. For any column vector $ x \in \R^d $ or $ k \in \N^d $, we use $ |\cdot| $ to denote the $ \ell_1 $-norm. Given any column vector $ x $ or matrix $ A $, we write $ x^T $ or $ A^T $ for their transposes. For any two column vectors $ x, y \in \R^d $, we define their inner product as
$$ \langle x, y \rangle = x^T y. $$
We say that a vector $ x = (x_r)_{r \in \S} \in \R^d $ satisfies $ x \leq y $ if and only if $ x_r \leq y_r $ for all $ r \in \S $. For any $ \delta \in \R $ and $ x \in \R^d $, we write $ x - \delta = (x_r - \delta)_{r \in \S} $. We use $ \| \kappa \|_\infty $ to denote the max norm of the matrix $ \kappa $.

%% file: cpp.tex
\section{The multi-dimensional compound Poisson process}\label{sec:cpp}
This section introduces a multi-dimensional compound Poisson process (MCPP), which will serve as a representation for the type-configurations of all connected components in the random graph $G_n$.

We first demonstrate that the distribution of the MCPP's jumps, conditional on its terminal state, is identical to the distribution of the type-configurations of all components in the random graph $G_n$. This key distributional equivalence allows us to establish the following:
\begin{itemize}
    \item[-] The largest jump in this conditional MCPP exhibits the same distributional behavior as the type $\cal{C}^n_{\rm max}$of the giant component in $G_n$.
    \item[-] The empirical measure of the conditional MCPP's jumps shares the same distribution as the empirical measure of component types, $\{t_n(k),k\in \N^d\}$.
    \item[-] The total number of jumps in the conditional MCPP is distributed identically to the total number of components, $C_n$.
\end{itemize}
Consequently, our main results can be concisely reformulated in terms of this conditional multi-type compound Poisson process.

\subsection{Definitions of  the multi-dimensional compound Poisson process}
\begin{defi}\label{defipk}
For each $k \in \mathbb{N}^d$, consider a random graph with $d$ vertex types whose type set is $\mathcal{S}$. 
The graph contains $k_r$ vertices of type-$r$ for every $r \in \mathcal{S}$, and each pair consisting of a type-$r$ vertex and a type-$s$ vertex is connected independently with probability $\kappa(r,s)/n$. 
Let $p_n(k)$ denote the connection probability of this graph.
\end{defi}
\begin{defi}\label{process}
For all $n$ fixed, let $\{X^n_i\}_{i\ge 1}$ be a sequence of $i.i.d$ random vectors in $$\bS_{n}:=\{k:k_r\in\{1,\cdots,\lfloor\mu^n_r n\rfloor\},r\in\S\}$$ with  common law
\begin{equation}\label{define}
    \P(X^{n}=k)=
    \frac{1}{Z_n}
    \frac{n^{|k|-1}p_n(k)\prod_{s\in\S}\left({\mu^n_s}\prod_{r\in\S}\left(1-\frac{\kappa(r,s)}{n}\right)^{\mu^n_rn-\frac{k_r}{2}}\right)^{k_s}}
    {\prod_{s\in\S}k_s!},
\end{equation}
where
\[
    Z_n=\sum_{k\in\bS_{n}}
    \frac{n^{|k|-1}p_n(k)\prod_{s\in\S}\left({\mu^n_s}\prod_{r\in\S}\left(1-\frac{\kappa(r,s)}{n}\right)^{\mu^n_rn-\frac{k_r}{2}}\right)^{k_s}}
    {\prod_{s\in\S}k_s!}
\]
is the normalizing constant. Let $N(Z_nn)$ be a Poisson process with intensity $Z_n n$ and independent of the sequence $\{X^{n}_i,i\in\N\}$. 
\end{defi}

 \subsection{Connection to the multi-type Random Graph}
    As established in Andreis et al.~\cite{LDPinhomo}, the distribution of the empirical measure $\{t_n(k), k\in\N^d\}$ is given by
      \begin{multline}\label{tnpb}
        \P(t_n(k)=\gamma_k, \forall k\in\bS_{n})\\ =\left(\prod_{r\in\S}\left(n\mu^n_r\right)!\right)\times\prod_{k\in\N^d}\frac{1}{\gamma_k!}\left(\frac{p_n(k)}{\prod_{r\in\S}k_r!}\prod_{s\in\S}\left(\prod_{r\in\S}\left(1-\frac{\kappa(r,s)}{n}\right)^{\mu^n_rn-\frac{k_r}{2}}\right)^{k_s}\right)^{\gamma_k},
        \end{multline}
        where $\gamma_k\in\N $ must satisfy the mass constraint $\sum_{k\in\N^d}k\gamma_k=\mu^nn$.
        
Hence, we can show in the following proposition that
the distribution of the empirical measure $\{t_n(k)\}$ in the graph $G_n$ is equal to the distribution of the empirical measure of the jumps in the compound Poisson process $\{X^{n}_1,X^{n}_2,\dots,X^{n}_{N(Z_nn)}\}$ conditioned on its terminal state $\left\{\sum_{i=1}^{N(Z_n n)}X_i^{n}=\mu^n n\right\}$. For the largest jump under the $\ell^1-$norm in the compound Poisson process $\{X_i^{n}\}$, with some abuse of notation, we denote its type configuration by $\max_{1\leq i\leq N(nZ_n)}X_i^n$.
\begin{prop}\label{connection}
    For all $(\gamma_k)\in\N^{\bS_n}$ such that $\sum_{k\in \bS_n}k\gamma_k=\mu^n n$, 
    \begin{multline*}
        \P\left(t_n(k)=\gamma_k, \forall k\in\bS_{n}\right)
=\P\left(\sum_{i=1}^{N(nZ_n)}\ind{X_i^{n}=k}=\gamma_k,\forall k\in\bS_n\bigg|\sum_{i=1}^{N(nZ_n)}X_i^{n}=\mu^n n\right).
    \end{multline*}
    Moreover,
    \begin{multline}\label{Pn}
         \mathbb{P}\left(\sum_{i=1}^{N(nZ_n)}X_i^n=\mu^n n\right)=
         \frac{e^{-Z_n n}
         n^n
         \prod_{r,s}{(1-\frac{\kappa(r,s)}{n})^{\frac{n^2\mu^n_s\mu^n_r}{2}}\prod_s{(\mu^n_s)^{n\mu^n_s}}}}
         {\prod_{s\in\mathcal{S}}(n\mu^n_s)!}.
    \end{multline}
    For any $k\in \bS_n$,
    \[\P\left(\mathcal{C}^n_{\max}=k\right)=\P\left(\max_{1\leq i\leq N(nZ_n)}X_i^n=k\bigg|\sum_{i=1}^{N(nZ_n)}X_i^{n}=\mu^n n\right).
    \]
    For all $l\in\N^+$,
    \[
    \P\left(C_n=l\right)=\P\left(N(nZ_n)=l\bigg|\sum_{i=1}^{N(nZ_n)}X_i^n=\mu^n n\right).
    \]
\end{prop}
\begin{proof}
  By direct calculation of Poisson random variables, we have
    \begin{equation*}
     \begin{aligned} &\P\left(\sum_{i=1}^{N(nZ_n)}\ind{X_i^n=k}=\gamma_k,k\in\bS_n\bigg|\sum_{i=1}^{N(Z_n n)}X_i^n=\mu^n n\right) \\
     =&\frac{1}{\mathbb{P}\left(\sum_{i=1}^{N(nZ_n)}X_i^{n}=\mu^n n\right)}
     \P\left(N(Z_n n)=\sum_k\gamma_k\right)
     \frac{\left(\sum_k \gamma_k\right)!}{\prod_k \gamma_k!}
     \prod_k(\P(X^{n}=k))^{\gamma_k}\\
     =&
     \frac
     {   e^{-Z_nn}
         n^n
         \prod_{r,s}{\left(1-\frac{\kappa(r,s)}{n}\right)^{\frac{n^2\mu^n_s\mu^n_r}{2}}\prod_s{(\mu^n_s)^{n\mu^n_s}}}
     }
     {    \prod_{s\in\mathcal{S}}(n\mu^n_s)!
     \mathbb{P}\left(\sum_{i=1}^{N(nZ_n)}X_i^n=\mu^n n\right)
     }
     \P(t_n(k)=\gamma_k,\forall k\in\mathbb{S}_{n}).
     \end{aligned}
    \end{equation*}
    Comparing this result with the expression for $\P(t_n(k)=\gamma_k,\forall k\in\mathbb{S}_{n})$ in \eqref{tnpb},  the first two results are obtained.
 Once we know the distribution of empirical measure, we know the distribution of $\cal{C}^n_{\rm max}$ and $C_n$ immediately.
\end{proof}

\subsection{Re-statement of our main results via the conditional MCPP}
We now reformulate our main results (Theorems~\ref{thm1} through~\ref{thm5}) using the newly established conditional multi-type compound Poisson process framework.
\begin{theorem}\label{mainthm1}
    Suppose $\Sigma(\kappa, \mu)> 1$, then for any sequence $(a_n)$ such that $n^{1/4}\ll a_n\ll \sqrt{n}$ and $|\mu^n-\mu|=O(\frac{1}n)$, under the conditioning event 
\[
\left\{\sum_{i=1}^{N(Z_n n)}X_i^n=\mu^n n\right\},
\]
the sequence of random vectors
\[
\frac{1}{a_n\sqrt{n}}\left(\max_{1\leq i\leq N(nZ_n)}X_i^n-(\mu-\cc)n\right),
\]
satisfies the moderate deviation principle in $\R^{d}$ with speed $a_n^2$ and rate function $\cal{I}(x)$, where $\cal{I}(x)$ is given in Theorem~\ref{thm1}.
\end{theorem}

\begin{theorem}\label{em}
 Suppose $\Sigma(\kappa, \mu)> 1$, then for any fixed $k\in\N^d$, any sequence $(a_n)$ such that $n^{1/4}\ll a_n\ll \sqrt{n}$ and $|\mu^n-\mu|=O(\frac{1}n)$, under the conditioning event  
\[
\left\{\sum_{i=1}^{N(Z_n n)}X_i^n=\mu^n n\right\},
\]
the sequence of random variables
\[
\left(
\frac{\sum_{i=1}^{N(Z_n n)}\mathbbm{1}_{\{X_i^n=k\}}-h(k)n}
{a_n\sqrt{n}}
\right)
\]
satisfies the moderate deviation principle in $\R$ with speed $a_n^2$ and rate function $J_k(x)$, where $J_k(x)$ is given in Theorem~\ref{thm2}.
\end{theorem}

\begin{theorem}\label{total}
Suppose $\Sigma(\kappa, \mu)> 1$, then for any sequence $(a_n)$ such that $n^{1/4}\ll a_n\ll \sqrt{n}$ and $|\mu^n-\mu|=O(\frac{1}n)$, under the conditioning event  
\[
\left\{\sum_{i=1}^{N(Z_n n)}X_i^n=\mu^n n\right\},
\]
the sequence of random variables
\[
\left(
\frac{N(Z_nn)-\left(|\cc|-\frac{1}{2}\langle\cc,\kappa\cc\rangle\right)n}
{a_n\sqrt{n}}
\right)
\]
satisfies the moderate deviation principle in $\R$ with speed $a_n^2$ and rate function $i(x)$,  where $i(x)$ is given in Theorem~\ref{thm3}.

\end{theorem}

Theorems~\ref{mainthm1},~\ref{em} and~\ref{total} will be proven in Section \ref{sec:super}. The proofs for the following theorems, concerning the subcritical case, are omitted as they closely parallel the proofs of Theorems~\ref{em} and~\ref{total}.

\begin{theorem}\label{sub1}
Suppose $\Sigma(\kappa, \mu)< 1$, then for any fixed $k\in\N^d$, any sequence $(a_n)$ such that $\sqrt{\log n}\ll a_n\ll \sqrt{n}$ and $|\mu^n-\mu|=O(\frac{1}n)$, under the conditioning event 
\[
\left\{\sum_{i=1}^{N(Z_n n)}X_i^n=\mu^n n\right\},
\]
for all $k\in \N^d$, the sequence of random variables
\[
\left(
\frac{\sum_{i=1}^{N(Z_n n)}\mathbbm{1}_{\{X_i^n=k\}}-h(k)n}
{a_n\sqrt{n}}
\right)
\]
satisfies the moderate deviation principle with speed $a_n^2$ and rate function
$\tilde{J}_k(x)$, where $\tilde{J}_k(x)$ is given in Theorem~\ref{thm4}.
\end{theorem}

\begin{theorem}\label{sub2}
   Suppose $\Sigma(\kappa, \mu)< 1$, then for any sequence $(a_n)$ such that $\sqrt{\log n}\ll a_n\ll \sqrt{n}$ and $|\mu^n-\mu|=O(\frac{1}n)$, under the conditioning event 
\[
\left\{\sum_{i=1}^{N(Z_n n)}X_i^n=\mu^n n\right\},
\]
 the sequence of random variables
\[
\left(
\frac{N(Z_nn)-\left(1-\frac{1}{2}\langle\mu,\kappa\mu\rangle\right)n}
{a_n\sqrt{n}}
\right)
\]
satisfies the moderate deviation principle with speed $a_n^2$ and rate function
$\tilde{i}(x)$,  where $\tilde{i}(x)$ is given in Theorem~\ref{thm5}.
\end{theorem}

%% file: connect.tex
\section{Asymptotic behaviors of the compound Poisson process}\label{sec:connect}

In this section, we analyze the asymptotic behaviors of the compound Poisson process $\{X^{n}_1,X^{n}_2,\dots,X^{n}_{N(Z_nn)}\}$ introduced in Definition~\ref{process}. We begin by recalling some fundamental properties of the connection probability $p_n(k)$ (Definition~\ref{defipk}).  Subsequently,, we examine the limiting distribution of a truncated jump, and finally, we estimate the terminal condition
\begin{equation}\label{tercon}
\P\left(\sum_{i=1}^{N(Z_nn)}X_i^{n}=\mu^n n\right),
\end{equation}
by comparing it with the tail distribution of the non-truncated jump  $X^{n}$.

\subsection{Bounds for the connection probability.}
We cite two lemmas bounds for $p_n(k)$ established by Andreis et al.~\cite{LDPinhomo}.
\begin{lemma}[Lemma~3.5~and~Lemma~4.8 in~\cite{LDPinhomo}]\label{est}
    For any $n\in\N$ and $k\in \N^d$,
    \begin{equation}\label{estimate}
    \left(1-\frac{\Vert\kappa\Vert_\infty}{n}\right)^{\frac{|k|^2}2}\leq\frac{p_n(k)}{n^{1-|k|}\tau(k)}\leq 1,
    \end{equation}
    where $\tau(k)$ is defined in~\eqref{tau(k)}.
\end{lemma}

\begin{lemma}[Lemma~3.4 in~\cite{LDPinhomo}]\label{lm3}
    Fix $k\in\N^d$,  as $n\rightarrow\infty$,
    \begin{equation}\label{p-tau}
    p_n(k)=n^{-(|k|-1)}\tau(k)(1+o(1)).
    \end{equation}
\end{lemma}

We now establish sharper bounds for the connection probability $p_n(k)$. These refined estimates are essential for analyzing the terminal condition~\eqref{tercon} in Proposition~\ref{asymp}. The proofs, which adapt the arguments found in~\cite{LDPinhomo}, are provided in Appendix~\ref{appen} to maintain the flow of the main text.

\begin{lemma}[Estimation of the connection probabilities(upper bound)]\label{esti2p}
For any sequence $(k^n)\subset\N^d$ with $k_s^n=O(n)$ for each $s\in\S$, we have
\[
p_n(k^n)\leq \prod_{s\in\S}\left(\sqrt{2\pi k^n_s}\right){\exp\left(\frac{1}{2n}\sum_{s\in\mathcal{S}}(\kappa k^n)_s\right)}
{\prod_{s\in\mathcal{S}}(1-e^{-\frac{1}{n}(\kappa k^n)_s})^{k^n_s}} \left(1+O\left(\frac{1}{n}\right)\right).
\]
\end{lemma}

\begin{lemma}[Estimation of the connection probabilities(lower bound)]\label{p_lb}
In the supercritical regime $\Sigma(\kappa,\mu)>1$, for any sequence $(a_n)$ satisfying $n^{\frac{1}{4}}\ll a_n\ll \sqrt{n}$, any $M>0$ and any $\beta\in\R^d$, define $$k^n(\beta)=n(\mu-c)+a_n\sqrt{n}\beta+O(1).$$
Then, we have
\[
        \liminf_{n\to\infty}\inf_{|\beta|\leq M}\frac{1}{a_n^2}\log \prod_{s\in\S}\left(\frac{n}{\sqrt{k^n_s(\beta)}}\right)\exp\left(-\langle k^n(\beta),\log (1-e^{-\frac{1}{n}(\kappa k^n(\beta) )}\rangle\right)p_n(k^n(\beta))\geq 0.
\]

\end{lemma}

\subsection{Asymptotic behaviors of the jumps}
\label{sec:jumps}

From the estimate~\eqref{p-tau}, for any fixed $k\in\N^d$, we have
\begin{equation}\label{esti1p}
\lim_{n\rightarrow\infty}Z_n\P(X^{n}=k)=h(k),
\end{equation}
where the distribution of $X^{n}$ is defined in \eqref{define}.
If $Z_n$ converges to a finite limit, then the random vectors $X^{n}$ converge in law to a limiting random vector $X$ on $\N^d$, such that for all $k\in\N^d$, the probability mass function $\P(X=k)$ is proportional to $h(k)$.
Therefore, before establishing the asymptotic distribution of $X^{n}$, we need the following basic properties of the functions $h(k)$ for $k\in\N^d$ and the limiting behavior of $Z_n$. Specifically, we will show that 
   \begin{equation}\label{defx}
    \P(X=k)=\frac{h(k)}{|c|-\frac{1}{2}c^T\kappa c}.
    \end{equation}
We note that in the classical Erd\H{o}s R\'enyi random graph ($d=1$), $\{kh(k)/(\sum_{\ell} \ell h(\ell))\}_{k\in\N}$ corresponds to  a Borel distribution on $\N$. For  properties of the Borel distribution, we refer to the paper~\cite{MR0111078}.
\begin{lemma}\label{cor}\hfill
\begin{itemize}
\item[-] When $\Sigma(\kappa,\mu)>1$, let $c$ be the solution to the characteristic function~\eqref{cha}. 
\item[-] When $\Sigma(\kappa,\mu)<1$, set $c=\mu$. 
\end{itemize}
Then $h(k)$ is given by
\begin{equation}\label{defhk}
    h(k)=\tau(k)\prod_{s\in\S}\frac{(e^{-(\kappa \cc)_s}\cc_s)^{k_s}}{k_s!}.
\end{equation}
Furthermore, the following properities hold for $\{h(k)\}$:
\begin{enumerate}[label=(\roman*)]
    \item Total mass:
    \begin{equation}\label{cor1}
  \sum_{k\in\N^d}h(k)=|\cc|-\frac{1}{2}\cc^T\kappa \cc,
\end{equation}
\item First moments: for any $r\in \cal{S}$,
\begin{equation}\label{cor2}
\sum_{k\in\N^d}k_rh(k)=\cc_r.
\end{equation}
\item Second moments: let $\Phi$ be a symmetric matrix with entries $\Phi_{rs}=\sum_{k\in\mathbb{N}^d}k_rk_sh(k)$ for all $r,s\in\cal{S} $. Then $\Phi$ is positive definite, and its inverse is
\begin{equation}\label{cor3}
  \Phi^{-1}=(D_{\cc})^{-1}-\kappa.
\end{equation}
\item High-order moments: for any $\gamma\in\N^d$, the $\gamma$-moment is finite:
\begin{equation}\label{highmoment}
    \sum_{k\in\N^d}h(k)\prod_{s\in\S}k_s^{\gamma_s}<\infty.
\end{equation}

\end{enumerate}

\end{lemma}
\begin{proof}
    From Proposition 4.2 in \cite{LDPinhomo}, we know that, for any fixed $r\in\S$, 
$\{\P_r(Y=k)=\frac{k_r}{\cc_r}h(k)\}$ defines a probability measure on $\N^d$, which immediately yields the equation~\eqref{cor2}. 

For the second moments, let  $g_s(c)=e^{-(\kappa c)_s}c_s$ for any $s\in\S$.
   From \eqref{cor2}, we have:
    \[
\cc_r=\sum_{k\in\N^d}\frac{\tau(k)}{\prod_{s\in\S}k_s!}k_r\prod_{s\in\S}g_s(\cc)^{k_s}.
\]
Combining with~\eqref{defhk}, we have 
\[
\Phi_{rs}=\sum_{k\in\N^d}\frac{\tau(k)}{\prod_{s\in\S}k_s!}k_rk_s\prod_{s\in\S}g_s(\cc)^{k_s}=g_s(\cc)\frac{\partial \cc_r}{\partial g_s}(\cc).
\]
In matrix form, it is
\[
\Phi=\big(\frac{\partial c}{\partial g}\big)D_{g(c)},
\]
where $g(c)=(g_s(c),s\in\cal{S})$. The Jacobian matrix of $g$ with respect to $c$ is:
\[
\left(\frac{\partial g}{\partial c}\right)=D_{g(c)}D^{-1}_c\left(I-D_c\kappa\right).
\]
Since $\Sigma(\kappa,c)<1$,  the largest eigenvalue of $D_c\kappa$ is less than $1$, implying that
$\frac{\partial g}{\partial c}$ is invertible. Applying the inverse function theorem, we find that the matrix $\Phi$ is also invertible, with the desired result
\[
\Phi^{-1}=-\kappa +D_{\cc}^{-1}.
\]

For the total mass~\eqref{cor1}, the $d=1$ case is known. It is a fundamental property of the Borel distribution proved by Aldous and Pitman \cite{aldous1998tree}. We proceed by induction on $d$.  For any $0\le x\le \cc_d$, define the vector
\[
\tilde{c}(x):=(\cc_1,\cc_2,\dots,\cc_{d-1},x)^T,
\]
and the function
\[
\tilde{h}_{x}(k)=\tau(k)\prod_{s\in\cal{S}}\frac{(e^{-(\kappa \tilde{c}(x))_s}\tilde{c}(x)_s)^{k_s}}{k_s!},
\]
with convention $0^0=1$. Note that $h(k)=\tilde{h}_{c_d}(k)$ for all $k$. By direct calculation,  for all $x>0$,
\[
\frac{\diff }{\diff x}\tilde{h}_x(k)=\left(-(\kappa k)_d+\frac{k_d}{x}\right)\tilde{h}_x(k).
\]
By using~\eqref{cor2}, we have
\begin{equation}\label{sum1}
\sum_{k}(\kappa k)_d\tilde{h}_x(k)=\sum_{s=1}^{d-1}\kappa(d,s)\cc_s+\kappa(d,d)x,
\end{equation}
and
\begin{equation}\label{sum2}
\sum_{k}\frac{k_d}{x}\tilde{h}_x(k)=1.
\end{equation}
Therefore, by solving the differential equation, we get:
\[
\sum_{k}\tilde{h}_{\cc_d}(k)=\sum_{k}\tilde{h}_{0}(k)-\sum_{s=1}^{d-1}\kappa(d,s)\cc_s\cc_d-\frac{1}{2}\kappa(d,d)(\cc_d)^2+\cc_d.
\]
Noticing that $\tilde{h}_0(k)\neq 0$ only for  $k$ satisfying $k_d=0$. Let
\[
k'=(k_1,\dots,k_{d-1})^T
\]
and $\kappa'$ be the $(d-1)\times (d-1)$ symmetric matrix with $\kappa'(i,j)=\kappa(i,j)$ for all $1\le i,j\le d-1$. Then the term $\sum_k\tilde{h}_0(k)$ simplifies to a $(d-1)$-dimensional sum,
\[
\sum_{k}\tilde{h}_0(k)=\sum_{k'}\tau'(k')\prod_{i=1}^{d-1}\frac{(e^{-\sum_{j=1}^{d-1}\kappa'(i,j)\cc_j}\cc_i)^{k_i'}}{k_i'!},
\]
which, by the inductive hypothesis, equals 
\[
\sum_{i=1}^{d-1}|\cc_i|-\frac{1}{2}\sum_{i,j=1}^{d-1}\cc_i\kappa(i,j)\cc_j.
\]
Combining these steps results in the desired identity \eqref{cor1}.

The existence of high-order moments \eqref{highmoment} follows from the moment generating function $G_r(t)=\E \left[\prod_{s\in\S}t_s^{Y_s}\right]$, where the measure $P_r$ is defined above.  For any $t\in[0,\infty)^d$,
\[
G_r(t)=\sum_{k\in\N^d}\frac{k_r}{\cc_r}h(k)\prod_{s\in\S}t_s^{k_s}=\sum_{k\in\N^d}\frac{k_r}{\cc_r}\tau(k)\prod_{s\in\S}\frac{\left(e^{-(\kappa c)_s}c_st_s\right)^{k_s}}{k_s!}.
\]
Lemma 4.5 in \cite{LDPinhomo} shows that $G_r(t)$ is analytic when $t_s=1$ for all $s\in\cal{S}$. Therefore, for $\gamma\in\N^d$ such that $\gamma_r\geq 1$, $\sum_{k\in\N^d}h(k)\prod_{s\in\S}k_s^{\gamma_s}<\infty$. Since the choice of $r$ is arbitrary, the finiteness of the moments is established.
\end{proof}

\begin{remark}
  [Relation with multi-type Poisson Galton Watson tree]
Recall that we have defined a multi-type Poisson Galton-Watson process with parameters $\kappa$ and $\mu$ in Section~\ref{branching}. Its type-configuration  $Y=(Y_i,i\in\cal{S})$ of the total population obeys the distribution
\[
\P_r(Y=k)=\frac{k_r}{\cc_r}h(k).
\]
It is straightforward to verify that: 
$$\E_i Y_j=\sum_k\frac{k_jk_i}{\cc_i}h(k)=\frac{1}{\cc_i}\Phi_{ij}.$$
Therefore, the matrix $M:=D_{\cc}^{-1}\Phi $ serves as the expected offspring matrix, \emph{i.e.}, the average matrix of the process.
By the fundamental property of the multi-type branching process, we have 
\[
M =I+\kappa D_{\cc} M,
\]
which readily yields equation~\eqref{cor3}.

The quantity~\eqref{cor1} also admits an interpretation within the multi-type Poisson Galton-Watson process.  In the classical one-type case ($d=1$), it has been shown by Pitman and Aldous~\cite{aldous1998tree} that this quantity~\eqref{cor1} corresponds to the expectation of reciprocal of the  total progeny in the Galton-Watson tree.
For the multi-type case, we can see from the definition of $h(k)$ and $Y$ that
    \[
\sum h(k)=\sum_s \cc_s\E_s \frac{1}{|Y|}.
    \]
This indicates that  if the multi-type Poisson Galton-Watson tree starts with the initial distribution $c/|c|=(c_s/|c|,s\in\cal{S})$, then the expectation of reciprocal of the total progeny in the multi-type Galton-Watson tree should be 
\[
\E\frac{1}{|Y|}=1-\frac{1}{2|\cc|}{\cc}^T\kappa \cc.
\]
We believe that this result can be proved by pruning the multi-type Galton-Watson tree  in a manner similar to the approach used by Aldous and Pitman~\cite{aldous1998tree}.

\end{remark}

We now analyze the asymptotic behaviors of the jump vector $X^n$ (Definition \ref{process}) as $n\to\infty$.  For technical reasons, we introduce a truncated jump $X^{n,\alpha}$, defined for any $\alpha\in\R_+^d$ with $0<\alpha_r<\mu_r^n$, and for all $k\in\bS_{\alpha n}:=\{k:k_r\in\{1,\cdots,\lfloor \alpha_r n\rfloor\},r\in\S\}$:
   \[
   \P(X^{n,\alpha}=k)=\frac{Z_n\P(X^{n}=k)}{Z_n^{\alpha}},
   \]
   and $ \P(X^{n,\alpha}=k)=0$ for all $k\notin\bS_{\alpha n}$, where
   \begin{equation}\label{zalp}
    Z_n^\alpha=\sum_{k\in\bS_{\alpha n}}Z_n\P(X^{n}=k).
    \end{equation}
   In the following of the paper, when $\Sigma(\kappa,\mu^n)>1$, we  use  $c^n$ to denote the solution to the characteristic equation with respect to  $\mu^n$,
   \begin{equation}\label{chan}
        c^n_ie^{-(\kappa c^n)_i}=\mu_i e^{-(\kappa\mu^n)_i}, \forall i\in\S,
   \end{equation}
   with $0\le c^n_i\le \mu^n_i$ for all $i\in\S$ and $\Sigma(\kappa,c^n)\le 1$. On the other hand, when $\Sigma(\kappa,\mu^n)\leq1$, we have $c^n=\mu^n$. We will use the following estimate from  Lemma 4.5 in \cite{LDPinhomo} in the proof.
\begin{lemma}\label{4.5inLDP}
    For any $r\in\S$ and arbitrary $\nu\in(0,\infty)^d$, as $m\rightarrow\infty$, 
    \[
\sum_{k\in\N^d,|k|=m}\tau(k)k_r\prod_{s\in\S}\frac{(\nu_se^{-(\kappa\nu)})^{k_s}}{k_s!}=e^{o(m)}e^{-m\left(\Sigma(\kappa,\nu)-\log\Sigma(\kappa,\nu)-1\right)}.
    \]
\end{lemma}

\begin{prop}\label{moment}
When $\Sigma(\kappa,\mu)\neq 1$, under the Assumption~\ref{assump1}, and assuming $|\mu^n-\mu|=O(\frac1n)$, for all $0<\alpha\leq\mu^n$ and any  $r\in\S$, we have
\begin{enumerate}[label=(\roman*)]
    \item Normalization Constant:
    \begin{equation}\label{zn}
      Z_n^{\alpha}=|\cc|-\frac{1}{2}\cc^T\kappa \cc+O\left(\frac1{n}\right).
    \end{equation}
    \item First Moments:
     \[
    \E X^{n,\alpha}_r=\frac{\cc_r}{|\cc|-\frac{1}{2}\cc^T\kappa \cc}+O\left(\frac1{n}\right),
\]
    \item Second Moments: let $\Psi^\alpha$ be the $d\times d$ symmetric matrix with entries $\Psi^\alpha_{ij}=\E X^{n,\alpha}_iX^{n,\alpha}_j$ for all $i,j\in\S$. Then we have 
\begin{equation}\label{Psi}
\Psi^\alpha_{ij}=\frac{1}{|\cc|-\frac{1}{2}\cc^T\kappa \cc}\Phi_{ij}+O\left(\frac{1}{n}\right),
\end{equation}
where $\Phi$ is the matrix given by~\eqref{cor3}.
\item Probability Mass Function: for any $k\in\N^d$,  
    \begin{equation}\label{X=k}
    \P(X^{n,\alpha}=k)=\frac{h(k)}{|\cc|-\frac{1}{2}\cc^T\kappa \cc}+O\left(\frac1{n}\right).
    \end{equation}
    \item Logarithmic Moment Generating Function: there exists $\eta\in(0,\infty)^d$ such that
\[
\sup_n\log\E \exp\langle\eta, X^{n,\alpha}\rangle<\infty.
\]
\end{enumerate}
\end{prop}
\begin{proof}
From the definition \eqref{zalp} of $Z_n^\alpha$ and the connection probability estimate in \eqref{estimate}, we have the following upper and lower bounds:
    \begin{multline}
  \label{znalpha}
       \sum_{k\leq\alpha n}\bigg(\frac{\tau(k)}
    {\prod_{s\in\S}k_s!}\left(1-\frac{\|\kappa\|_\infty}{n}\right)^{\frac{|k|^2}{2}}\prod_{s\in\S}\left(\prod_{r\in\S}\left(1-\frac{\kappa(r,s)}{n}\right)^{\mu^n_rn-\frac{k_r}{2}}\mu^n_s\right)^{k_s}\bigg)\\
        \leq Z^\alpha_n\leq
        \sum_{k\leq\alpha n}\bigg(\frac{\tau(k)}
        {\prod_{s\in\S}k_s!}\prod_{s\in\S}\left(\prod_{r\in\S}\left(1-\frac{\kappa(r,s)}{n}\right)^{\mu^n_rn-\frac{k_r}{2}}\mu^n_s\right)^{k_s}\bigg).
  \end{multline}
    For every $m\in\N_+$, let $f_n(m)$ be the sum in the upper bound restricted to $|k|=m$:
    \[
    f_n(m)=\sum_{k:|k|=m, k\leq\mu^nn }\frac{\tau(k)}
        {\prod_{s\in\S}k_s!}\prod_{s\in\S}\left(\prod_{r\in\S}\left(1-\frac{\kappa(r,s)}{n}\right)^{\mu^n_rn-\frac{k_r}{2}}\mu^n_s\right)^{k_s}.
    \]
We can bound $f_n(m)$ by
    \begin{equation}\label{leq00}
    \left(f_n(m)\right)^{\frac1{m}}   \leq\left(\sum_{k:|k|=m} h^n(k)\right)^{\frac{1}{m}}
    \left( \sup_{k:|k|=m}\prod_{s\in\S}\prod_{r\in\S}\left(1-\frac{\kappa(r,s)}{n}\right)^{-\frac{k_rk_s}{2|k|}}\right),
    \end{equation}
    where 
  \[
 h^n(k)=\tau(k)\prod_{s\in\S}\frac{(e^{-(\kappa\mu^n)_s}\mu^n_s)^{k_s}}{k_s!}.
   \]
 Since $|\mu^n-\mu|=O\left(\frac{1}{n}\right)$ and the function $\Sigma(\kappa,\mu)$ is continuous in $\mu$, there exists a $\tilde{\mu}\in\R_+^d$ such that $\Sigma(\kappa,\tilde{\mu})-\log\Sigma(\kappa,\tilde{\mu})-\frac12[\kappa]>1$ and  $e^{-(\kappa\Tilde{\mu})_s}\Tilde{\mu}_s>e^{-(\kappa\mu^n)_s}\mu^n_s,\ \forall s\in\S$ for  $n$ sufficiently large.
   As $n\rightarrow\infty$, the superior limit of the second term in~\eqref{leq00} can be controlled by the definition of $[\kappa]$ in Assumption~\ref{assump1},
    \[    \limsup_{n\rightarrow\infty}\sup_{k:|k|=m}\prod_{s\in\S}\prod_{r\in\S}\left(1-\frac{\kappa(r,s)}{n}\right)^{-\frac{k_rk_s}{2|k|}}\leq \exp\left(\frac{1}{2}\left\langle\frac{k}{|k|},\kappa\frac{k}{|k|}\right\rangle\right)\leq e^{\frac12[\kappa]}.
    \]
Thus, for  $n$ sufficiently large, for all $m\in\N$, we have
\[
\left(f_n(m)\right)^{\frac1m}\leq(f(m))^{\frac1m}:=\left(\sum_{k:|k|=m}\tilde{h}(k)\right)^{\frac{1}{m}}e^{\frac{1}{2}[\kappa]},
    \]
    where 
 \[
   \tilde{h}(k)=\tau(k)\prod_{s\in\S}\frac{(e^{-(\kappa\Tilde{\mu})_s}\Tilde{\mu}_s)^{k_s}}{k_s!}.
   \]    
Therefore, the upper bound in~\eqref{znalpha} gives us
\begin{equation}\label{znupp}
Z_n^{\alpha}\le \sum_{0\leq m\leq|\alpha n|}f_n(m)\le \sum_{0\leq m\leq|\alpha n|}f(m).    
\end{equation}
Now we use Lemma \ref{4.5inLDP} to obtain
    \[    f(m)=e^{\frac{1}{2}[\kappa]m}\sum_{k:|k|=m}\tau(k)\prod_{s\in\S}\left(\frac{(e^{-(\kappa\tilde{\mu})_s}\tilde{\mu}_s)^{k_s}}{k_s!}\right)\leq e^{o(m)}e^{-m\left(\Sigma(\kappa,\tilde{\mu})-\log\Sigma(\kappa,\tilde{\mu})-1-\frac{1}{2}[\kappa]\right)}.
    \]
Hence, under the Assumption~\ref{assump1}, we have
    \begin{equation}\label{leq1}
\limsup_{m\rightarrow\infty}\left(f(m)\right)^{\frac1{m}}<1,
    \end{equation}
which  implies $Z_n^\alpha<\infty$.

 Now we calculate the limit of $Z_n^\alpha$.    Let $g_n(k)$ denote the ratio of the $k$-th term in the upper bound~\eqref{znalpha} to $h^n(k)$:
 $$g_n(k):=\prod_{s\in\S}\left(\prod_{r\in\S}\left(1-\frac{\kappa(r,s)}{n}\right)^{\mu^n_rn-\frac{k_r}{2}}e^{\left(\kappa \mu^n\right)_s}\right)^{k_s}.$$
By using Taylor expansion,  we have
    \begin{multline}\label{definegn}
    \log g_n(k)\\
    =\frac{1}{n}\left(-\sum_{s\in\S}\sum_{r\in\S}\frac{\kappa (r,s)^2\mu^n_sk_r}{2}+\frac{k^T\kappa k}{2}\right)+\frac{1}{n^2}\left(\frac{\kappa(r,s)^2}{2}k_rk_s\right)+O\left(\frac1n\right).
    \end{multline}
Hence, for sufficiently large $n$, we have
    \begin{equation}\label{prop21}
    \begin{aligned}
    &\left|\sum_{k\leq\alpha n}\left(\frac{\tau(k)}
        {\prod_{s\in\S}k_s!}\prod_{s\in\S}\left(\prod_{r\in\S}\left(1-\frac{\kappa(r,s)}{n}\right)^{\mu^n_rn-\frac{k_r}{2}}\mu^n_s\right)^{k_s}-h^n(k)\right)\right|\\
        \leq&\sum_{k\leq\alpha n}{h^n}(k)\left|g_n(k)-1\right|\\
        \leq&\sum_{|k|\leq\sqrt{n}}\tilde{h}(k)\left|g_n(k)-1-\log g_n(k)\right|+\sum_{|k|\leq\sqrt{n}}\tilde{h}(k)\left|
        \log g_n(k)\right|\\
&\qquad\qquad\qquad\qquad\qquad+C_1\sum_{m>\sqrt{n}}f(m),
    \end{aligned}
    \end{equation}
 using the fact $\tilde{h}(k)>h^n(k)$ for all $k\in\N^d$.
    For any $|k|\le\sqrt{n}$,  we have $\log g_n(k)<\|\kappa\|_\infty /2+1$ for $n$ sufficiently large. In this case, there exists a constant $C_2>0$, such that the inequality~\eqref{prop21} is upper bounded by 
    \[
C_2\sum_{|k|\leq\sqrt{n}}\tilde{h}(k)\left|
        \log g_n(k)\right|+C_1\sum_{m>\sqrt{n}}f(m).
    \]
 Using the Lemma~\ref{cor}, we know that $\sum\Tilde{h}(k)$, $\sum k_r\Tilde{h}(k)$ and $\sum k_rk_s\Tilde{h}(k)$ are finite.    Hence, using \eqref{leq1} and the Taylor expansion~\eqref{definegn}, the equation~\eqref{prop21} is $O\left(\frac{1}{n}\right)$. 
    Similarly, for the lower bound in~\eqref{znalpha},
    \[    
    \sum_{k\leq\alpha n}\bigg|\frac{\tau(k)}
    {\prod_{s\in\S}k_s!}\left(1-\frac{\|\kappa\|_\infty}{n}\right)^{\frac{|k|^2}{2}}\prod_{s\in\S}\left(\prod_{r\in\S}\left(1-\frac{\kappa(r,s)}{n}\right)^{\mu^n_rn-\frac{k_r}{2}}\mu^n_s\right)^{k_s}-h^n(k)\bigg|
        =O\left(\frac1n\right).
    \]
Thanks to the Lemma~\ref{cor}, we conclude that
    \[
    Z_n^\alpha=\sum_{k\le \alpha n}h^n(k)+O\left(\frac1n\right)=|\cc^n|-\frac{1}{2}(\cc^n)^T\kappa \cc^n+O\left(\frac1n\right).
    \]
    Since both $(\mu,c)$ and $(\mu^n,c^n)$ satisfy the dual equation~\eqref{cha}, we know that $c^n-c=O(\frac{1}{n})$ by the assumption $\mu^n-\mu=O(\frac{1}{n})$. Thus,
    \[
    Z_n^\alpha=|\cc|-\frac{1}{2}\cc^{T}\kappa \cc+O\left(\frac1{n}\right).
    \]
Using the relation \eqref{highmoment}, we obtain the following results by a similar approach 
\[
    \E X^{n,\alpha}_r=\sum_{k\leq\alpha n}k_rh(k)=\frac{\cc_r}{|\cc|-\frac{1}{2}\cc^T\kappa \cc}+O\left(\frac1{n}\right).
\]
    and
\[
    \Psi^{n,\alpha}_{ij}=\sum_{k\leq\alpha n}k_rk_sh(k)=\frac{1}{|\cc|-\frac{1}{2}\cc^T\kappa \cc}\Phi_{ij}+O\left(\frac{1}{n}\right).
\]
    It remains to prove that the logarithmic moment generating function is finite. We  can follow a similar approach to  the upper bounds \eqref{znalpha} and \eqref{leq1} and show that there exist $\eta\in(0,+\infty)^d$ small enough and $m_0\in\N$ large enough such that 
\begin{multline*}
     Z^\alpha_n\E \exp\langle\eta,X^{n,\alpha}\rangle\le\sum_{m=1}^{\infty}\sum_{|k|=m}\exp\<\eta,k\>Z_n^{\alpha}\P(X^{n,\alpha}=k)\\
    \le\sum_{m=1}^{m_0}\sum_{k:|k|=m}\exp\<\eta,k\>Z_n^{\alpha}\P(X^{n,\alpha}=k)\\
    +
    \sum_{m=m_0}^{\infty}e^{m\|\eta\|_\infty}
 e^{-m\left(\Sigma(\kappa,\tilde{\mu})-\log\Sigma(\kappa,\tilde{\mu})-1-\frac{1}{2}[\kappa]\right)}<\infty.
\end{multline*}
\end{proof}

\subsection{Estimates of the terminal condition}
We now provide an estimation for the probability of the terminal condition,
\[
\P\left(\sum_{i=1}^{N(Z_nn)}X_i^{n}=\mu^n n\right)
\]
by comparison with the tail distribution of the jump vector  $X^{n}$ defined in definition~\ref{process}.

\begin{prop}[Asymptotic Behavior of the Terminal Condition]\label{asymp}
\hfill
\begin{enumerate}[label=(\roman*)]
    \item Subcritical Case $(\Sigma(\kappa,\mu)<1)$: for any sequence $a_n$ such that $a_n\gg (\log n)^\frac12$, we have
\[
\lim_{n\rightarrow\infty}\mdp \left(\sumn X_i^{n}=\mu^n n\right)=0.
\]
\item Supercritical Case $(\Sigma(\kappa,\mu)>1)$: for any sequence $n^{\frac{1}{4}}\ll a_n\ll\sqrt{n}$, any $\mu^n$ with $|\mu^n_r-\mu_r|=O(\frac 1n)$, $\forall r$ and any $M>0$, we have 
\begin{multline*}
  \lim_{n\rightarrow\infty}\sup_{|\beta|\leq M}\bigg|\frac{1}{a_n^2}\log \frac
{Z_nn\mathbb{P}(X^n=\lfloor n(\mu^n-\cc^n)+\beta a_n\sqrt{n}\rfloor)}
{\mathbb{P}(\sum_{i=1}^{N(nZ_n)}X_i^n=\mu^n n)}\\
+
\frac12\left\langle D_{\mu-\cc}^{-1}(I-D_{\cc} \kappa)\beta,(I-D_{\mu} \kappa)\beta\right\rangle
\bigg|=0.
\end{multline*}
\end{enumerate}

\end{prop}
\begin{proof}
 
In the subcritical case $\Sigma(\kappa,\mmu)<1$, we have $\cc=\mmu$. Using the formula~\eqref{Pn}, Stirling's approximation, and the estimation~\eqref{zn}, we obtain
\begin{multline*}
\lim_{n\to\infty}\mdp \left(\sumn X_i^{n}=\mu^n n\right)\\
=\lim_{n\to\infty}\frac{1}{a_n^2}\left(-Z_nn+|\cc^n|n-\frac{1}{2}n(\cc^n)^{T}\kappa\cc^n-\sum_{s\in\S}\log(\sqrt{2\pi n\mu^n_s})\right)
=0,
\end{multline*}
 when $a_n\gg (\log n)^{1/2}$.

In the supercritical case $\Sigma(\kappa,\mmu)>1$, we first use the relation \eqref{Pn} and Stirling's approximation to get, for all $k_s\gg 1$,~$\forall s\in\cal{S}$,
\begin{equation}\label{fra}
\begin{aligned}
&\log\frac{Z_n n \P(X^n=k)}{\P(\sum_{i=1}^{N(nZ_n)}X_i^n=\mmu^n n)}\\
=&\log\left(\prod_{r,s\in\cal{S}}\left(1-\frac{\kappa(r,s)}{n}\right)^{k_s\mu^n_rn-\frac{k_sk_r}{2}-\frac{n^2\mu^n_s\mu^n_r}{2}}\right)\\
     &+\log \left(p_n(k)\prod_{s\in\mathcal{S}}\left(
\frac{
n\mu^n_s
}
{k_s}\right)^{k_s}
e^{Z_nn+|k|-n|\mu^n|}\right)
+
\sum_{s\in\mathcal{S}}\log\left(\sqrt{
\frac{n\mu^n_s}
{k_s}
}
\right)+O(1).
\end{aligned}
\end{equation}

In particular, for $k^n(\beta)=\lfloor n(\mu^n-\cc^n)+\beta a_n\sqrt{n}\rfloor$, $\beta \in\R^d$, thanks to the fact $\sum_r\mu^n_r=1$, the first term in the right-hand side of the equation~\eqref{fra} becomes
\begin{equation}\label{39}
\begin{aligned}
&\log\left(\prod_{r,s\in\cal{S}}\left(1-\frac{\kappa(r,s)}{n}\right)^{k^n_s(\beta)\mu^n_rn-\frac{k^n_s(\beta)k^n_r(\beta)}{2}-\frac{n^2\mu^n_s\mu^n_r}{2}}\right)\\
&=\sum_{r,s}\frac{1}{2}\left(-\cc^n_r\cc^n_sn^2+\beta_s\cc^n_ra_nn^{\frac{3}{2}}+\beta_r\cc^n_sa_nn^{\frac{3}{2}}-\beta_r\beta_sa_n^2n\right)\log\left(1-\frac{\kappa(r,s)}{n}\right)\\
&=\frac{1}{2}(\cc^n)^T\kappa \cc^nn-\beta^T\kappa \cc^na_n\sqrt{n}+\frac{1}{2}\beta^T\kappa\beta a_n^2+O(1).
\end{aligned}
\end{equation}
 On the other hand, by applying the Lemma~\ref{esti2p}, the second term in  equation~\eqref{fra} has an upper bound,
\begin{multline}\label{23}
\log \left(p_n(k^n(\beta))\prod_{s\in\mathcal{S}}\left(
\frac{
n\mu^n_s
}
{k^n_s(\beta)}\right)^{k^n_s(\beta)}
e^{Z_nn+|k^n(\beta)|-n|\mu^n|}\right)
\\
\le-\frac{1}{2}\cc^T\kappa \cc n+|\beta|a_n\sqrt{n}+\sum_{s\in\mathcal{S}}k^n_s(\beta)(\log\frac{n}{k^n_s(\beta)})+\sum_{s\in\mathcal{S}}k^n_s(\beta)\log \mu^n_s\\
+\frac{d}{2}\log (2\pi n)+\frac{1}{2n}\sum_{s\in\mathcal{S}}(\kappa k^n(\beta))_s
+\sum_{s\in\mathcal{S}}k^n_s(\beta)\log(1-e^{-\frac{1}{n}(\kappa k^n(\beta))_s})+O(1).
\end{multline}
Moreover, using the dual equation~\eqref{chan} and Taylor's expansion, we obtain
\begin{equation}\label{41}
\begin{aligned}
&\sum_{s\in\mathcal{S}}k^n_s(\beta) \log\left(\frac{n}{k^n_s(\beta)}\left(1-e^{-\frac{1}{n}(\kappa k^n(\beta))_s}\right)\right)\\
=&\sum_{s\in\mathcal{S}}k^n_s(\beta) \log\left(\frac{1}{\mu^n_s-\cc^n_s+\beta_s\frac{a_n}{\sqrt{n}}}\left(1-\frac{\cc^n_s}{\mu^n_s}e^{-(\kappa \beta)_s\frac{a_n}{\sqrt{n}}}\right)\right)\\
=&
-\sum_{s\in\mathcal{S}}k^n_s(\beta)\log \mu^n_s+
(-|\beta|+(\cc^n)^T\kappa\beta)a_n\sqrt{n}\\
&-\sum_{s\in\mathcal{S}}\left\{
\frac{1}{2}\cc^n_s(\kappa \beta)_s^2
+\frac{1}{2}\frac{(\beta_s-\cc^n_s(\kappa\beta)_s)^2}{\mu^n_s-\cc^n_s}
\right\}a_n^2
+O_{\beta}\left(1+\frac{a_n^3}{\sqrt{n}}\right).
\end{aligned}
\end{equation}
In conclusion, we get
\begin{multline*}
\frac{1}{a_n^2}\log \frac{Z_nn\mathbb{P}(X^n=\lfloor n(\mu^n-\cc^n)+\beta a_n\sqrt{n})\rfloor}{\mathbb{P}(\sum_{i=1}^{N(nZ_n)}X_i^n=\mu^n n)}\\
\leq-\frac12\left\langle D_{\mu^n-\cc^n}^{-1}(I- D_{\cc^n}\kappa)\beta,(I-D_{\mu^n} \kappa)\beta)\right\rangle
+O_{\beta}\left(\frac{1}{a_n^2}+\frac{a_n}{\sqrt{n}}\right)\\
=-\frac12\left\langle D_{\mu-\cc}^{-1}(I- D_{\cc}\kappa)\beta,(I-D_{\mu} \kappa)\beta)\right\rangle
+O_{\beta}\left(\frac{1}{a_n^2}+\frac{a_n}{\sqrt{n}}\right).
\end{multline*}
For the lower bound, we combine the lower bound in Lemma~\ref{p_lb} with the terms in \eqref{fra},~\eqref{39} and~\eqref{41} to deduce the lower bound
\begin{multline*}
\liminf_{n\to\infty}\inf_{|\beta|\leq M}\frac{1}{a_n^2}\log\frac{Z_n n \P(X^n=\lfloor n(\mu^n-\cc^n)+\beta a_n\sqrt{n})\rfloor)}{\P(\sum_{i=1}^{N(nZ_n)}X_i^n=\mmu^n n)}
\\\geq\liminf_{n\to\infty}\inf_{|\beta|\leq M}\frac{1}{a_n^2}\Bigg(
(\cc^n)^T\kappa \cc^nn-\beta^T\kappa \cc^na_n\sqrt{n}+\frac{1}{2}\beta^T\kappa\beta a_n^2
\\+\sum_{s\in\S}\log\left(\frac{\sqrt{k_s^n(\beta)}}{n}\right)+\sum_{s\in\S}k_s^n(\beta)\log\left(1-e^{-\left(\kappa\frac{k^n(\beta)}{n}\right)_s}\right)
\\+k_s^n(\beta)\log\left(\frac{n\mu^n_s}{k_s^n(\beta)}\right)+Z_nn+|k^n(\beta)|-n
\Bigg)\\
=-\frac12\left\langle D_{\mu-\cc}^{-1}(I- D_{\cc}\kappa)\beta,(I-D_{\mu} \kappa)\beta)\right\rangle
.
\end{multline*}
The proof is complete.
\end{proof}

%% file: mdpcpp.tex
\section{MDP for the compound Poisson process.}\label{sec:mdpcpp}
In this section, we first apply the G\"artner-Ellis theorem to show the moderate deviation principle for a general multi-dimensional compound Poisson process. We then give  explicit rate functions related to the moderate deviations of the truncated process $\{X_1^{n,\alpha},X_2^{n,\alpha},\dots,X_{N(Z_n n)}^{n,\alpha}\}$.
\begin{lemma}\label{generalMDP}
  Let $a_n$ be a sequence such that $1\ll a_n\ll\sqrt{n}$. Let $(\xi^{n,k})_{k\in\N}$ be a sequence of  $i.i.d.$ random vectors in  $\R^d$. Suppose that their moments satisfy:
  $$\E \xi^{n,1}=m+O(\frac{1}{\sqrt{n}}),$$
  and $$\E(\xi^{n,1}_i-m_i)(\xi^{n,1}_j-m_j)=\sigma_{ij}+O(\frac{1}{\sqrt{n}}),$$
  for some $m\in\R^d$ and some covariance matrix $C=(\sigma_{ij})_{d\times d}$. Let $b_n\in\R^d$ be a sequence such that $|b_{n}|=O\left(\sqrt{n}+a_n^2\right)$.
  Suppose that there exists $R>0$ such that for all $\eta\in\R^d$ and $|\eta|<R$, $\sup_n\log \E\exp(\eta\cdot \xi^n)<\infty$. 
\begin{enumerate}[label=(\roman*)]
    \item   Let $N(\lambda_n n)$ be an independent Poisson process with $\lambda_n=\lambda+O(\frac{1}{n})$. The sequence of random variables 
    $$\frac{\sum_{k=1}^{N(\lambda_n n)}\xi^{n,k}-\lambda m n+b_n}{a_n\sqrt{n}}$$
    satisfies the moderate deviation principle in $\R^d$ with speed $a_n^2$ and the rate function 
    \[
    \sup_{l\in\R^d}\left\{\langle l,x\rangle-\frac\lambda{2}\langle l, (C+mm^T)l\rangle\right\}.
    \]
    Moreover, when the matrix $C+mm^T$ is positive-definite, the rate function equals  
    $${\frac{1}{2\lambda}\<x,(C+mm^T)^{-1}x\>.}$$
    \item For any $\lambda>0$ and $u\in\R$, the sequence of random variables 
    $$\frac{\sum_{k=1}^{\lfloor\lambda n+ua_n\sqrt{n}\rfloor}\xi^{n,k}-(\lambda n+ua_n\sqrt{n})m+b_n}{a_n\sqrt{n}}$$
    satisfies the moderate deviation principle in $\R^d$ with speed $a_n^2$ and the rate function
    \[
    \sup_{l\in\R^d}\left\{\langle l,x\rangle-\frac\lambda{2}\langle l, C l\rangle\right\}.
    \]
    Moreover, when the matrix $C$ is positive-definite, the rate function equals
    $${\frac{1}{2\lambda}\<x,C^{-1}x\>.}$$
\end{enumerate}
\end{lemma}

\begin{proof}
     For any $z\in\R^d$, by Taylor expansion and  the dominated convergence theorem, we have  estimates of the scaled cumulant generating functions,
    \begin{equation}\label{mdppoi}
        \begin{aligned}
            &\frac{1}{a_n^2}\log\E\left(\exp\left\{\frac{a_n}{\sqrt{n}}\left<z, \sum_{k=1}^{N(\lambda_n n)}\xi^{n,k}-\lambda n m+b_n\right>\right\}\right)\\
            =&\frac{\lambda n+O(\sqrt{n})}{a_n^2}\left(
            \frac{a_n}{\sqrt{n}}\<z,\E\xi^n
            \>+\frac{1}{2}\left(\frac{a_n}{\sqrt{n}}\right)^2\E\left(\<z,\xi^n\>\right)^2
            +O\left(\left(\frac{a_n}{\sqrt{n}}\right)^3\right)\right)\\
            &\qquad\qquad\qquad-\lambda\frac{\sqrt{n}}{a_n} \<z,m\>+\frac{1}{a_n\sqrt{n}}\<z, b_n\>
            \\
            =&\frac{1}{2}\lambda \left(z^TCz+\<z, m\>^2\right)
            +O_z\left(\frac{1}{a_n}+\frac{a_n}{\sqrt{n}}\right),
        \end{aligned}
    \end{equation}
    and
    \begin{equation}\label{mdpsum}
        \begin{aligned}
            &\frac{1}{a_n^2}\log\E\left(\exp\left\{\frac{a_n}{\sqrt{n}}\left\langle z, \sum_{k=1}^{\lfloor\lambda n+ua_n\sqrt{n}\rfloor}\xi^{n,k}-(\lambda n +u a_n\sqrt{n})m+b_n\right\rangle\right\}\right)\\
            =&\frac{\lfloor\lambda n+ua_n\sqrt{n}\rfloor}{a_n^2}\left(
            \frac{a_n}{\sqrt{n}}\langle z,\E\xi^n\rangle
            +\frac{1}{2}\left(\frac{a_n}{\sqrt{n}}\right)^2 z^TCz
            +O\left(\left(\frac{a_n}{\sqrt{n}}\right)^3\right)\right)\\
            &\qquad\qquad\qquad
            -\frac{\lambda\sqrt{n}+ua_n}{a_n} \langle z, m\rangle+\frac{1}{a_n\sqrt{n}}\langle z,b_n\rangle 
            \\
            =&\frac{1}{2}\lambda z^TCz
            +O_z\left(\frac{1}{a_n}+\frac{a_n}{\sqrt{n}}\right).
        \end{aligned}
    \end{equation}
    Thanks to the G\"{a}rtner-Ellis theorem, we derive the two principles of moderate deviations. See Theorem 3.7.1 in Dembo and Zeitouni~\cite{dembo2009large} for a similar result. When the matrices are positive definite, we can apply the Legendre-Fenchel transform and obtain the corresponding explicit rate functions.
\end{proof}

We now establish the positive definiteness or invertibility of the matrices that appear in the rate functions  of the moderate deviations.
\begin{lemma}
  The matrix in Theorem~\ref{thm1},
$(I-\kappa D_c)D_{\nu}(I-D_c\kappa),$
  is  positive definite.
\end{lemma}
\begin{proof}
Since all the eigenvalues of $D_v$ are strictly positive and  $\Sigma(\kappa,c)<1$, the matrix $(I-\kappa D_c)D_{\nu}(I-D_c\kappa)$ is positive definite.  
\end{proof}
\begin{lemma}\label{lmbk}
      The matrix  $\Phi-kk^Th(k)$ in Theorem~\ref{thm2} and in Theorem~\ref{thm4} is positive definite.
\end{lemma}
\begin{proof}
Let $\tilde{X}:=X\ind{X\neq k}$, where $X$ is defined in~\eqref{defx}. Recalling Lemma~\ref{cor}, we have that the matrix,
\[\Phi-kk^Th(k)=\left(|c|-\frac{1}{2}c^T\kappa c\right)\E\left(\tilde{X} \tilde{X}^T\right),\]
is a positive semi-definite matrix. 
Assuming this matrix is singular,  there exists a non-zero eigenvector $a\in\R^d$, such that $(\Phi-kk^Th(k))a=0$. This implies $\E (a^T\tilde{X})^2=0.$ On the other hand, for this vector $a$, we can always find $k'\in\N^d\setminus\{0\}$ and $k'\neq k$, such that $|a^Tk'|> 0$. Therefore,
\[
\E (a^T\tilde{X})^2\ge (a^Tk')^2\P\left(\tilde{X}=k'\right)=(a^Tk')^2\frac{h(k')}{\sum_{k''\neq k}h(k'')}>0, 
\]
which leads to a contradiction. In conclusion, the matrix $(\Phi-kk^Th(k))$ must be positive definite.
\end{proof}
\begin{lemma}\label{lmb}
 The matrix $\Phi-\frac{\cc\cc^T}{|c|-\frac{1}{2}\cc^T\kappa \cc}$ in Theorem~\ref{thm3} and in Theorem~\ref{thm5} is positive definite.
\end{lemma}
\begin{proof}
 We first observe that the matrix
 \[
 \Phi-\frac{\cc\cc^T}{|c|-\frac{1}{2}\cc^T\kappa \cc}=(|c|-\frac{1}{2}\cc^T\kappa \cc)\textrm{Cov}(X)
 \]
is positive semi-definite. We only need to check that it is invertible.
The matrix $\Phi^{-1}\cc\cc^T$ has rank $1$ and then has one non-zero eigenvalue.
Therefore, we have
\[
\textrm{det} \left(I-\frac{\Phi^{-1} \cc \cc^T}{|\cc|-\frac12\<\cc,\kappa\cc\>}\right)=1-\textrm{tr} \left(\frac{\Phi^{-1} \cc \cc^T}{|\cc|-\frac12\<\cc,\kappa\cc\>}\right).
\]
Recalling that $\Phi^{-1}=D_{\cc}^{-1}-\kappa$, the above determinant becomes:
\[
    \frac{\frac{1}{2}\<\cc,\kappa\cc\>}{|\cc|-\frac12\<\cc,\kappa\cc\>},
  \]
  which is  strictly positive. Hence, the matrix $I-\frac{\Phi^{-1} \cc \cc^T}{|\cc|-\frac12\<\cc,\kappa\cc\>}$ is invertible. In conclusion, $\Phi-\frac{\cc\cc^T}{|c|-\frac{1}{2}\cc^T\kappa \cc}$ is positive definite.
\end{proof}
\begin{lemma}\label{pdefinite1}
    The matrix $A+B_k$ in Theorem \ref{thm2} is positive definite.
\end{lemma}

\begin{proof}
    Recall that the matrix $\Phi^{-1}=D_c^{-1}-\kappa$ is a positive definite matrix. Hence, we can use the Sherman-Morrison formula and obtain,
    \[
    B_k=\Phi^{-1}+h(k)\frac{\Phi^{-1} kk^T\Phi^{-1}}{1-h(k)k^T\Phi^{-1} k},
    \]
    and then
    \begin{multline*}
A+B_k=(I-\kappa D_\cc )D^{-1}_{\mu-\cc}(I-D_\mu \kappa)+B_k\\
=\Phi^{-1}D_cD^{-1}_{\mu-c}D_\mu\Phi^{-1}+h(k)\frac{\Phi^{-1} kk^T\Phi^{-1}}{1-h(k)k^T\Phi^{-1} k}.
    \end{multline*}
    It is clear that the matrix $\Phi^{-1}D_cD^{-1}_{\mu-c}D_\mu\Phi^{-1}$ is positive definite and $\Phi^{-1} kk^T\Phi^{-1}$ is positive semi-definite. 
By a straightforward calculation and the positive definiteness of the matrices $\Phi^{-1}$ and $B_k$, we find
 \[
   1-h(k)k^T\Phi^{-1} k =\frac{k^T\Phi^{-1}k}{k^TB_k k}>0.
    \]
    Therefore, the matrix $A+B_k$ is positive definite.
  
\end{proof}

\begin{lemma}\label{pdefinite2}
The matrix $A+B$ in Theorem~\ref{thm3} is positive definite.
\end{lemma}
\begin{proof}
By the Sherman-Morrison formula for $B$, we have
\[
A+B=\Phi^{-1}D_cD^{-1}_{\mu-c}D_\mu\Phi^{-1}+\frac{\Phi^{-1} cc^T\Phi^{-1}}{\frac 12\langle \cc, \kappa\cc\rangle}.
\]
For a reason similar to the proof of the Lemma \ref{pdefinite1}, it is positive definite.

\end{proof}

Now we give the following moderate deviation principles related to our compound Poisson process. These results will be used in the next section to prove our main results. We remark that under the condition $X_i^n\le \alpha^n n$, the law of $X_i^n$ is given by the truncated distribution $X_i^{n,\alpha}$.

\begin{prop}\label{MDPcppp}
     Under the Assumption~\ref{assump1}, for any sequence $a_n$ such that $1\ll a_n\ll\sqrt{n}$, suppose that $|\mu^n-\mu|=O\left(\frac{1}{n}\right)$ and $\alpha\in\R_+^d$ satisfy $\alpha_r\le\mu_r$ for all $r\in\S$. Let $\{\alpha^n\}\subset\R_+^d$ be a sequence such that $\alpha^n\rightarrow\alpha$ as $n\rightarrow\infty$ and $0<\alpha^n_r\le\mu^n_r$ for all $r\in\S$. Under condition $X_{i,r}^n\leq \alpha^n_r n$ for all $i\geq 1$ and $r\in\S$, 
    \begin{enumerate}
        \item [(i)]  the sequence of random vectors 
        $$\frac{1}{a_n\sqrt{n}}\left(\sum_{i=1}^{N(Z_n n)}X_i^{n}-\cc^n n\right)$$
        satisfies the moderate deviation principle in $\R^d$ with speed $a_n^2$ and the rate function  
        $$j_1(x)=\frac{1}{2}\<x, \Phi^{-1} x\>,$$
        where $\Phi^{-1}=D_{\cc}^{-1}-\kappa$;
        \item [(ii)] the sequence of random vectors 
        $$\frac{1}{a_n\sqrt{n}}\left(\sum_{i=1}^{\lfloor(|\cc^n|-\frac12\langle\cc^n,\kappa\cc^n\rangle)n+a_n\sqrt{n}x\rfloor}X_i^{n}-\cc^n n-\frac{\cc^n a_n\sqrt{n}x}{|\cc^n|-\frac 12\langle \cc^n, \kappa\cc^n\rangle}\right)$$
        satisfies the moderate deviation principle in $\R^d$ with speed $a_n^2$ and the rate function 
        $$
       j_2(x)= \frac12 \bigg\langle x,\left(\Phi-\frac{\cc\cc^T}{|\cc|-\frac12\langle\cc,\kappa\cc\rangle}\right)^{-1}x \bigg\rangle;
        $$
        
      \item [(iii)]the sequence of random vectors 
        \[\frac{1}{a_n\sqrt{n}}\left(\sum_{i=1}^{N(Z_n n)}X_i^{n}\mathbbm{1}_{\{X^{n}_i\neq k\}}-(\cc^n-kh(k)) n\right)
        \]
        satisfies the moderate deviation principle with speed $a_n^2$ and rate function  
        \[
        j_3(x)=\frac{1}{2}\<x,(\Phi-kk^Th(k))^{-1}x\>.
        \]

    \end{enumerate}
\end{prop}
\begin{proof}
  Since $\mu$ and $\cc$ satisfy the dual equation \eqref{cha}, we can regard $\mu$ as a function of $c$ and differentiate accordingly. Given that the derivative is bounded, it follows that $c^n_r-c_r=O(\frac{1}{n})$ when $\mu^n_r-\mu_r=O\left(\frac{1}{n}\right)$. 
  Under the condition $X_i^n\leq \alpha^n n$ for all $i\geq 1$, the law of $X_i^n$ is $X_i^{n,\alpha^n}$. Applying Lemma~\ref{generalMDP} and Proposition~\ref{moment}, the three moderate deviation principles hold with the following rate functions,
  \begin{align*}
(i)\qquad&j_1(x)=\sup_{\ell\in\R^{d}}\left\{\langle \ell,x \rangle-\frac{1}{2}\langle \ell, \Phi \ell\rangle\right\},\\ (ii)\qquad&j_2(x)=\sup_{\ell\in\R^{d}}\left\{\langle \ell,x \rangle-\frac{1}{2}\left\langle \ell, \left(\Phi-\frac{\cc\cc^T}{|\cc|-\frac12\langle\cc,\kappa\cc\rangle}\right) \ell\right\rangle\right\},\\    (iii)\qquad&j_3(x)=\sup_{\ell\in\R^{d}}\left\{\langle \ell,x \rangle-\frac{1}{2}\left\langle \ell, (\Phi-kk^Th(k))\ell\right\rangle\right\}.
  \end{align*}
  Thanks to Lemma~\ref{cor}, Lemma~\ref{lmb} and Lemma~\ref{lmbk}, the matrices $\Phi$, $\Phi-\frac{\cc\cc^T}{|\cc|-\frac12\langle\cc,\kappa\cc\rangle}$ and $\Phi-kk^Th(k)$ are all positive definite. Consequently, we obtain the desired results.

\end{proof}

The following technical lemma will also be used in the next section.
\begin{lemma}\label{term3}
    For any $\alpha\in \R^d$ such that $\alpha_r>0, \forall r \in\S$ and any $a_n\gg 1$, we have  
    \[
    \lim_{n\rightarrow\infty}\frac{1}{a_n^2}\log\P\bigg((X_{i}^n)_s\leq \alpha_s n,\ \forall s\in\S,\  1\leq \forall i\leq N(Z_n n)\bigg)=0.
    \]
\end{lemma}
\begin{proof}
    For any fixed $r\in\S$ , by a straightforward calculation and Proposition \ref{moment}, we have
   \begin{multline*}
        \left|\frac{1}{a_n^2}\log\P\bigg((X_{i}^n)_s\leq \alpha_s n,\ \forall s\in\S,\   1\leq \forall i\leq N(Z_n n)\bigg)\right|=\frac{Z_nn}{a_n^2} \P(X_{s}^n> \alpha_s n,\ \forall s\in\S )\\
        \leq \frac{Z_nn}{a_n^2} \P(X^n_r>\alpha_r n)\le
        \frac{Z_n\E X_r^n}{\alpha_r a_n^2}\to 0,\qquad \textrm{as }n\to\infty.
    \end{multline*}
    
\end{proof}

%% file: super.tex
\section{Proofs of the moderate deviation principle in the supercritical case}\label{sec:super}
In this section, we will prove the main results in the super critical case. For any vector $\beta \in \R^d$ and radius $\delta>0$, 
we denote by $B_\delta(\beta)$ and $\overline{B_\delta(\beta)}$ the $\ell_1$-open  and closed balls, respectively, centered at $\beta$.

\subsection{Proofs for the largest component}
We now present the proofs of the main result, Theorem~\ref{mainthm1}, which characterizes the moderate deviation behavior of the largest component size.



\begin{lemma}[upper bound]\label{upmax}\label{up1}
For any scaling sequence $n^{1/4}\ll a_n\ll \sqrt{n}$ and any $\beta\in\R^{d}$, the following conditional superior limit holds:
\[
\begin{aligned}
    &\lim_{\delta\to 0}\limsup_{n\rightarrow\infty}\frac{1}{a_n^2}\log\P\left(\frac{1}{a_n\sqrt{n}}\left(\max_{1\leq i\leq N(nZ_n)}X_i^n-(\mu-\cc)n\right)\in \overline{B_\delta(\beta)}\bigg|\sum_{i=1}^{N(nZ_n)}X_i^n=\mu^n n\right)\\
   \leq& -\frac{1}{2}\left\langle\left(I-D_\cc\kappa\right)\beta,D_v\left(I-D_\cc\kappa\right)\beta\right\rangle.
\end{aligned}
\]
\end{lemma}

\begin{proof}
Let $\tX^n$ be an independent copy of $X^n_i$. For any $x\in\R^d$ and $n$ sufficiently large, define $k^n(x)=\lfloor(\mu^n-\cc^n)n-xa_n\sqrt{n}\rfloor$. Since the sequence $\{X_i^n\}$ consists of independent and identically distributed random variables, we establish the following upper bound on the joint probability:
\begin{align*}
    &\qquad\P\left(\max_{1\leq i\leq N(nZ_n)}X_i^n=k^n(x),\sum_{i=1}^{N(nZ_n)}X_i^n=\mu^n n \right)\\
        &\leq\sum_{k=0}^\infty\frac{e^{-Z_n n}(Z_n n)^{k+1}}{(k+1)!}(k+1)\\
        &\qquad\qquad\times \P\left(\tX^n=k^n(x), \sum_{i=1}^kX^n_i=\mu^n n-k^n(x), |X^n_i|\leq |k^n(x)|, 1\leq i\leq k\right)\\
  &\qquad\qquad\qquad\qquad \leq Z_n n \P\left(\tX^n=k^n(x)\right)\P\left(\sum_{i=1}^{N(Z_nn)}X^n_i=\mu^n n-k^n(x)\right).
  \end{align*}
Consequently, for any $\delta>0$, the scaled  logarithm   of the conditional probability can be bounded as:
\begin{multline*}
        \frac{1}{a_n^2}\log \P\left(\frac{1}{a_n\sqrt{n}}\left(\max_{1\leq i\leq N(nZ_n)}X_i^n-(\mu^n-\cc^n)n\right)\in\overline{\B}\bigg|\sum_{i=1}^{N(Z_nn)}X^n_i=\mu^n n\right)\\
      \leq\frac{1}{a_n^2}\log \sup_{x\in\overline{\B}}\frac{Z_n n\P(X^n=k^n(x))}{\P\left(\sum_{i=1}^{N(nZ_n)}X_i^n=\mu^n n\right)}\\
      +\frac{1}{a_n^2}\log\P\left(\frac{1}{a_n\sqrt{n}}\left(\sum_{i=1}^{N(nZ_n)}X_i^n-\cc^n n\right)\in\overline{\B}\right).
    \end{multline*}
The limit of the first term on the right-hand side is determined by Proposition~\ref{asymp}. The superior limit of the second term is obtained from the moderate deviation principle  for the sequence $\sum_{i=1}^{N(nZ_n)}X_i^n$, as stated in Proposition~\ref{MDPcppp}. By combining these limits and utilizing the fact that $(\mu^n_r-c^n_r)-(\mu_r-c_r)=O\left(\frac{1}{n}\right)$, we derive the desired upper bound:
\[
    \begin{aligned}
      &\lim_{\delta\to 0}\mdpsup\left(\frac{1}{a_n\sqrt{n}}\left(\max_{1\leq i\leq N(nZ_n)}X_i^n-(\mu-\cc)n\right)\in\overline{\B}\bigg|\sum_{i=1}^{N(nZ_n)}X^n_i=\mu^n n\right)\\
      \leq& -\frac{1}{2}\left(\left\langle (I-D_\cc\kappa)\beta,D^{-1}_{\mu-\cc}(I-D_\mu\kappa)\beta\right\rangle+\left\langle (I-D_\cc \kappa)\beta,D_\cc^{-1}\beta\right\rangle\right)\\
        =& -\frac{1}{2}\left\langle\left(I-D_\cc\kappa\right)\beta,D_v\left(I-D_\cc\kappa\right)\beta\right\rangle.
    \end{aligned}
\]

\end{proof}

\begin{lemma}[lower bound]\label{low1}
    For any scaling sequence $n^{\frac{1}{4}}\ll a_n\ll\sqrt{n}$ and $\beta>0$,
    the corresponding conditional inferior limit satisfies the following:
\begin{equation*}
\begin{aligned}
    &\lim_{\delta\to 0}\mdpinf\left(\frac{1}{a_n\sqrt{n}}\left(\max_{1\leq i\leq N(nZ_n)}X_i^n-(\mu-\cc)n\right)\in B_\delta(\beta)\bigg|\sum_{i=1}^{N(nZ_n)}X_i^n=\mu^n n\right)\\
    \geq& -\frac{1}{2}\left\langle\left(I-D_\cc\kappa\right)\beta,D_v\left(I-D_\cc\kappa\right)\beta\right\rangle.
\end{aligned}
\end{equation*}
\end{lemma}

\begin{proof}
 Let $k^n(x)=\lfloor(\mu^n-\cc^n)n-xa_n\sqrt{n}\rfloor$. For any vector $\eps\in\R^d$ such that $0<\eps_r<\mu_r-\cc_r, \forall r\in\S$ and for $n$  sufficiently large, we utilize a technique analogous to the proof of the upper bound. However, here we restrict the event to ensure that one specific component achieves the maximum size, $k^n(x)$, while all other components are bounded away from this value by $\eps$. In this case, we establish the lower bound:
    \begin{equation*}
    \begin{aligned}
        &\P\left(\max_{1\leq i\leq N(nZ_n)}X_i^n=k^n(x),\sum_{i=1}^{N(nZ_n)}X_i^n=\mu^n n \right)\\
        \geq&Z_n n\P\left( X^n=k^n(x)\right)\\
        \times&\P\left(\sum_{i=1}^{N(Z_n n)}X_i^n=\lceil \cc^n n+a_n\sqrt{n}x\rceil, X_i^n\leq\left(\mu^n-\cc^n-\eps\right)n, 1\leq\forall i\leq N(Z_nn)\right).
    \end{aligned}
    \end{equation*}
Utilizing the conditional probability definition, we derive the following lower bound on the scaled logarithm of the conditional probability:
    \begin{equation*}
        \begin{aligned}
            &\mdp\left(\frac{1}{a_n\sqrt{n}}\left(\max_{1\leq i\leq N(nZ_n)}X_i^n-(\mu^n-\cc^n)n\right)\in\B\Bigg|\sum_{i=1}^{N(nZ_n)}X^n_i=\mu^n n\right)\\          
            \geq&\frac{1}{a_n^2}\Bigg\{\log\inf_{x\in\B}\frac{Z_n n\P(X^n=k^n(x))}{\P\left(\sum_{i=1}^{N(nZ_n)}X_i^n=\mu^n n\right)}\\            
            +&\log\P\left(\frac{1}{a_n\sqrt{n}}\left(\sum_{i=1}^{N(nZ_n)}X_i^n-\cc^n n\right)\in\B\bigg|X_i^n\leq(\mu^n-\cc^n-\eps)n, 1\leq i\leq N(Z_nn)\right)\\
            +&\log\P(X_i^n\leq(\mu^n-\cc^n-\eps)n, 1\leq i\leq N(Z_nn))
            \Bigg\}.
        \end{aligned}
    \end{equation*}
The limit of the first term is obtained through Proposition~\ref{asymp}. The inferior limit of the second term is derived from the conditional MDP result presented in Proposition~\ref{MDPcppp}. By Lemma~\ref{term3}, the third term is shown to be asymptotically vanishing as $n\rightarrow \infty$. 
Finally, combining these limits and incorporating the assumption that $(\mu^n_r-c^n_r)-(\mu_r-c_r)=O\left(\frac{1}{n}\right)$, we successfully establish the requisite lower bound.
\end{proof}

\begin{lemma}[Exponential Tightness]\label{et1}
    For any $M>0$ and $n\in\N^+$, define a compact set
    \[
    A_n(M):=\left\{\bigg|\max_{1\leq i\leq N(Z_n n)} X_i^n-(\mu^n-\cc^n)n \bigg|\le Ma_n\sqrt{n}\right\}.
    \]
    Then for any scaling sequence $\{a_n\}_{n\geq 1}$ such that $n^{\frac{1}{4}}\ll a_n\ll \sqrt{n}$, 
    \[
    \begin{aligned}
    \mdpsup\left(A_n(M)^c\Bigg|\sum_{i=1}^{N(Z_n n)}X_i^n=\mu^n n\right)&\\
    \leq -\frac{M^2\lambda_1}{2},
    \end{aligned}
    \]
    where $\lambda_1>0$ is the minimum eigenvalue of the matrix $(I-\kappa D_c)D_{\nu}(I-D_c\kappa)$.
\end{lemma}

\begin{proof}
    For all sequences $\{l_n\}_{n\geq 1}\subset\R^d$ such that $n\gg|l_n|>Ma_n\sqrt{n}$,  we apply the upper bound established in Lemma \ref{upmax}. This yields: 
    \[
    \begin{aligned}
        &\limsup_{n\rightarrow\infty}\frac{n}{l_n^2}
        \log\P
        \left(\max_{1\leq i\leq N(Z_n n)}X_i^n=\lfloor(\mu-\cc)n+l_n\rfloor\bigg|\sum_{i=1}^{N(Z_nn)}X_i^n=\mu^n n\right)\\
        \leq&-\frac{\lambda_1}{2}.
    \end{aligned}
    \]
   On the other hand, for any $y\in \R^d_+\setminus\{0\}$ and any sequence satisfying $\limn \frac{l_n}{n}=y$, by the large deviation results proved in Andreis~\cite{LDPinhomo}, there exist $\eta_y>0$ such that
    \[
    \limsup_{n\rightarrow\infty}\frac1n\log\P\left(\max_{1\leq i\leq N(Z_n n)}X_i^n=\lfloor(\mu-\cc)n+l_n\rfloor\bigg|\sumn X_i^n=\mu^n n\right)\leq -\eta_y.
    \]
  Combining these results, we conclude that for any sequence $\{l_n\}_{n\geq 1}\subset\R^d$ such that $|l_n|>Ma_n\sqrt{n}$ and  $0\leq  l_n+(\mu-\cc)n\leq \mu n$, the superior limit of scaled logarithm of the following conditional probability is  bounded by
    \[
    \begin{aligned}
    &\mdpsup \left(\max_{1\leq i\leq N(Z_n n)}X_i^n=\lfloor(\mu-\cc)n+l_n\rfloor\bigg|\sumn X_i^n=\mu^n n\right)\\
    &\leq -\frac{M^2\lambda_1}{2}.
    \end{aligned}
    \]
    Finally, noting that the cardinality of the set $A_n (M )^c$ under the condition $\{\sum_{i=1}^{N(Z_nn)} X_i^n=\mu^n n\}$ is $O(n)$, and given the scaling assumption $n^{\frac{1}{2}}\ll a_n^2\ll n$, we apply the Laplace's method to derive the exponential tightness.
\end{proof}

\begin{proof}[Proof of Theorem~\ref{mainthm1}]
    Thanks to the Theorem~4.1.11 in Dembo and Zeitouni~\cite{dembo2009large}, by using the upper bound in Lemma~\ref{up1}, the lower bound in Lemma~\ref{low1} and the exponential tightness in Lemma~\ref{et1},  we have the  moderate deviation principle for all sequence $\{a_n\}$ satisfying $n^{1/4}\ll a_n\ll \sqrt{n}$. 
\end{proof}


\subsection{Proofs for the empirical measure.}

This section presents the proof of Theorem~\ref{em}, which establishes the moderate deviation principle for the empirical measure $\sum_{i=1}^{N(Z_nn)}\ind{X_i^n=k}$ for all $k\in\N^d$ conditioned on the total sum
$\sum_{i=1}^{N(Z_nn)}X_i^n=\mu^nn$.
\begin{proof}
The proof leverages the properties of independent Poisson Point Processes. Specifically, for any fixed vector $k\in\N^d$, the random variables $\sumn \mathbbm{1}_{\{X_i^n=k\}}$ (the count of random vectors equal to $k$) and $\sumn X_i^n\mathbbm{1}_{\{X_i^n\neq k\}}$ (the count of random vectors not equal to $k$) are independent. 
The conditional probability of the count $\sumn \mathbbm{1}_{\{X_i^n=k\}}$ being equal to $j$ for some  $j\in\N^+$ is given by
\begin{multline}\label{emp1}
    \P\left(\sumn \mathbbm{1}_{\{X_i^n=k\}}=j\bigg|\sumn X_i^n=\mu^nn\right)\\
=e^{-Z_n\P(X^n=k)n}
\frac{(Z_n\P(X^n=k)n)^j}{j!}
\frac{\P\left(\sumn X_i^n\mathbbm{1}_{\{X_i^n\neq k\}}=\mu^n n-j k\right)}
{\P\left(\sumn X_i^n=\mu^n n\right)}.
\end{multline}
Similar to the methodology employed in the proof for the largest component, we introduce the following event $A_n'(M)$ for any $M\in\mathbb{R}^+$ concerning the maximum value of the non-$k$ variables:
\[
A_n'(M):=\left\{
\left|
\max_{1\leq i\leq N(Z_nn)}X_i^n\mathbbm{1}_{\{X_i\neq k\}}-(\mu^n-\cc^n)n
\right|\leq M a_n\sqrt{n}
\right\}.
\]
By applying Lemma~\ref{et1}, we can readily demonstrate that
\begin{equation}\label{a'}
\lim_{M\to\infty}\limsup_{n\to\infty}\frac{1}{a_n^2}\log\P\left(A_n'(M)^c\Bigg|\sumn X_i^n=\mu^nn\right)=-\infty.
\end{equation}
For any measurable set $G\subset \R$, let $G_k=\{xk\in \R^d:x\in G\}$. Thanks to relation~\eqref{emp1}, the following probability
\[
\P\left(
\frac{1}{a_n\sqrt{n}}\left(
\sumn \mathbbm{1}_{\{X_i^n=k\}}-h(k)n
\right)\in G
\bigg| \sumn X_i^n=\mu^nn\right)
\]
is upper bound by the following two terms
\begin{multline}\label{empu1}
\P\left(
\frac{1}{a_n\sqrt{n}}\left(
\sumn \mathbbm{1}_{\{X_i^n=k\}}-h(k)n
\right)\in G
\right)\\
\times
\frac{\P\left(\frac{1}{a_n\sqrt{n}}\left(\sumn X_i^n\mathbbm{1}_{\left\{X_i^n\neq k\right\}}-n(\mu^n-kh(k))\right)\in G_k, A_n'(M)\right)}{\P\left(\sumn X_i^n=\mu^nn\right)}
\\
+\P\left(A_n'(M)^c\bigg|\sumn X_i^n=\mu^n n\right).
\end{multline}
Thanks to Chebyshev's inequality, inequality~\eqref{mdppoi} and Proposition~\ref{asymp}, for $y\in\R^d$ with $|y|< M$, $k^n(y)=\lfloor(\mu^n- \cc^n)n+y a_n\sqrt{n}\rfloor\in\N^d$ and any $z\in\R^d$, we have the estimate:
\begin{multline}\label{empupper}
    \frac{1}{a_n^2}\log\P\Bigg(\frac{1}{a_n\sqrt{n}}\left(\sumn X_i^n\mathbbm{1}_{\left\{X_i^n\neq k\right\}}-n(\mu^n-kh(k))\right)\in G_k,\\
   \max_{1\leq i\leq N(Z_nn)}X_i^n\mathbbm{1}_{\{X_i\neq k\}}=k^n(y)\Bigg)-\frac{1}{a_n^2}\log\P\left(\sumn X_i^n=\mu^nn\right)\\
    \leq\frac{1}{a_n^2}\log \frac{Z_nn \P\left(
    X^n=k^n(y)
    \right)}{\P\left(\sumn X_i^n=\mu^nn\right)}\hfill\\
    +\frac{1}{a_n^2}\log \P\left(\frac{1}{a_n\sqrt{n}}\left(\sumn X_i^n\mathbbm{1}_{\left\{X_i^n\neq k\right\}}-n(\cc^n-kh(k))\right)+y\in G_k\right)\\
    \leq
    \left(
    -
    \frac{1}{2}\langle
    D_{\mu-c}^{-1}(I-D_{c} \kappa)y,(I-D_{\mu}\kappa)y
    \rangle
    +\frac{1}{2}\langle z,(\Phi-kk^Th(k))z\rangle
    -
    \inf_{x\in G}\langle z,xk-y\rangle
    \right)\\
    +O_{z}\left(\frac{1}{a_n^2}+\frac{a_n}{\sqrt{n}}\right)+O_{M}\left(\frac{1}{a_n^2}+\frac{a_n}{\sqrt{n}}\right).
\end{multline}
\paragraph{\bf Upper bound} For a given $\beta\in\R$, we select the measurable set $G=\overline {B_\delta(\beta)}$ in the inequality~\eqref{empu1}. The first probability in \eqref{empu1} is estimated using the MDP results for the Poisson random variables, yielding the term $-\beta^2/h(k)$ in the rate function. 

In order to estimate the ratio term in~\eqref{empu1}, we notice that the cardinality of the event $A_n'(M)$ is $O\left((Ma_n\sqrt{n})^d\right)$ and $\log (a_n\sqrt{n})^d/a_n^2\to 0.$ Hence, we can choose the optimizing parameter $z_y=(\Phi-kk^T h(k))^{-1}(\beta k-y)$ for any given $y$ in~\eqref{empupper} and then give an upper bound for the ratio term:
\begin{multline}\label{empu2}
    \sup_{y:|y|\leq M}\left(
    -
    \frac{1}{2}\langle
    D_{\mu-c}^{-1}(I-D_{c} \kappa)y,(I-D_{\mu}\kappa)y
    \rangle
    +\frac{1}{2}\langle z_y,(\Phi-kk^Th(k))z_y\rangle
    -\langle z_y,\beta k-y\rangle
    \right)\\
    +C_k\delta M+O_M\left(\frac{\log( a_n\sqrt{n})}{a_n^2}\right)+O_M\left(\frac{1}{a_n^2}+\frac{a_n}{\sqrt{n}}\right)+O_{M}\left(\frac{1}{a_n^2}+\frac{a_n}{\sqrt{n}}\right),
\end{multline}
where $C_k$ is a positive constant  not depending on $M,\delta$ or $n$. By taking $M$ sufficiently large and invoking Lemma~\ref{pdefinite1}, which confirms the positive definiteness of $A+B_k$, we find that the supremum over $|y|\le M$ in~\eqref{empu2} is attained at  $$y=\beta(A+B_k)^{-1}B_kk$$ and equal to 
\[-\frac{\beta^2}{2}\left(k^T A(A+B_k)^{-1}B_kk\right),\]
where $A=(I-D_{c} \kappa)D_{\mu-c}^{-1}(I-D_{\mu} \kappa)$ and $B_k=(\Phi-kk^Th(k))^{-1}$.

Assembling these results together with~\eqref{a'} and taking the limits $n\to\infty$ and $\delta\to 0$, we conclude the moderate deviation upper bound:
\begin{equation}\label{needinver}
    \begin{aligned}
        &\lim_{\delta\rightarrow 0}\mdpsup 
        \left(
        \frac{1}{a_n\sqrt{n}}\left(
        \sumn \mathbbm{1}_{\{X_i^n=k\}}-h(k)n\right)\in \overline{\B}\bigg|\sumn X_i^n=\mu^nn
        \right)\\
        \leq&
        -\frac{\beta^2}{2}\left(\frac{1}{h(k)}+k^T A(A+B_k)^{-1}B_kk\right).
    \end{aligned}
\end{equation}

\paragraph{\bf Exponential tightness} 
By employing inequalities~\eqref{empu1} and~\eqref{empupper}, we establish the following bound for any $b>0$ and $M>0$,
\begin{multline*}
     \P\left(
    \frac{1}{a_n\sqrt{n}}\bigg|
    \sumn \mathbbm{1}_{\{X_i^n=k\}}-h(k)n\bigg|>b\Bigg|\sumn X_i^n=\mu^nn
    \right)\\
    \le  \P\left(
    \frac{1}{a_n\sqrt{n}}\bigg|
    \sumn \mathbbm{1}_{\{X_i^n=k\}}-h(k)n\bigg|>b
    \right)\sup_{|y|\le M}\frac{Z_nn \P\left(
    X^n=k^n(y)
    \right)}{\P\left(\sumn X_i^n=\mu^nn\right)}\\
+\P\left(A_n'(M)^c\bigg|\sumn X_i^n=\mu^n n\right).
\end{multline*}
The desired exponential tightness is then achieved by combining this bound with Proposition~\ref{asymp}, the limit established in~\eqref{a'}, and the basic properties of Poisson random variables. Specifically, we first choose $M$ sufficiently large, and subsequently take the limit as $b\to\infty$.

\paragraph{\bf Lower bound} For any $\beta\in\R$, we select the specific vector $y=\beta(A+B_k)^{-1}B_kk$ that was found to maximize the quadratic form in the upper bound analysis. For any $x \in\B$ and $0<\alpha<\mmu-\cc$,
\begin{multline*}
    \P\left(\sumn \mathbbm{1}_{\{X_i^n=k\}}=\lfloor h(k)n+xa_n\sqrt{n}\rfloor\bigg|\sumn X_i^n=\mu^nn\right)\\
    \ge \P\left(\left\{\sumn \mathbbm{1}_{\{X_i^n=k\}}=\lfloor h(k)n+xa_n\sqrt{n}\rfloor\right\}\bigcap D^n(\alpha)\bigg|\sumn X_i^n=\mu^nn\right),
\end{multline*}
where $D^n(\alpha)$ is the event
\[
D^n(\alpha):=\left\{\max_{1\le i\le N(Z_nn)}X_i^n=k^n(y), i^*=\arg\max_{1\le i\le N(Z_nn)}X_i^n,X_i^n\le \alpha n,\forall i\neq i^*\right\}.
\]
Similar to the proof of the largest component, this probability is lower bounded by 
\begin{multline*}
   \P\left(\sumn \mathbbm{1}_{\{X_i^n=k\}}=\lfloor h(k)n+xa_n\sqrt{n}\rfloor\right)\times\left( \frac{Z_nn\P(X^n=k^n(y))}{\P\left(\sum_{i=1}^{N(Z_nn)}X_i^n=\mu^n n\right)}\right)\\
    \times\P\left( \sumn X^n_i\mathbbm{1}_{\left\{X_i^n\neq k\right\}}=-knh(k)+c^n n-a_n\sqrt{n}(kx-y), X_i^n\leq n\alpha\right).
\end{multline*}
Applying Stirling's approximation to the Poisson random variable, we have
\[
\begin{aligned}
&\lim_{\delta\rightarrow 0}\liminf_{n\rightarrow\infty}\inf_{x\in\B}\frac{1}{a_n^2}\log\P\left(\sumn \mathbbm{1}_{\{X_i^n=k\}}=\lfloor h(k)n+xa_n\sqrt{n}\rfloor\right)\\
=&\lim_{\delta\rightarrow 0}\liminf_{n\rightarrow\infty}\inf_{x\in\B}\frac{1}{a_n^2}\log\left(e^{-Z_n\P(X^n=k)n}
\frac{(Z_n\P(X^n=k)n)^{\lfloor h(k)n+xa_n\sqrt{n}\rfloor}}{\lfloor h(k)n+xa_n\sqrt{n}\rfloor!}\right)\\
\geq& -\frac{\beta^2}{2h(k)}.
\end{aligned}
\]
Using Proposition \ref{MDPcppp} and Lemma \ref{term3}, we know
\[
\begin{aligned}
\lim_{\delta\rightarrow 0}\liminf_{n\rightarrow\infty}\frac{1}{a_n^2}\log\P\Bigg(\frac{1}{a_n\sqrt{n}}&\left( \sumn X_i^n\mathbbm{1}_{\{X_i^n\neq k\}}-n(\cc^n-kh(k))\right)\in B_{k\delta}(y-k\beta),\\
&X_i^n\leq\alpha n, \forall i\Bigg)\geq
-\frac{1}{2}\<y-k\beta,(\Phi-kk^T h(k))^{-1}(y-k\beta)\>.
\end{aligned}
\]
 Combining the terms above and Proposition \ref{asymp}, we establish the lower bound
\begin{multline}\label{totallow}
\lim_{\delta\rightarrow 0}\mdpinf 
        \left(
        \frac{1}{a_n\sqrt{n}}\left(
        \sumn \mathbbm{1}_{\{X_i^n=k\}}-h(k)n\right)\in \B\bigg|\sumn X_i^n=\mu^nn
        \right)\\
        \geq-\frac{\beta^2}{2h(k)}-  
        \frac{1}{2}\bigg(\langle
        D_{\mu-\cc}^{-1}(I-D_\cc \kappa)y,(I-D_\mu \kappa)y
        \rangle
        +
        \langle y-\beta k, (\Phi-kk^T h(k))^{-1}(y-\beta k)\rangle
        \bigg)\\
        =-\frac{\beta^2}{2}\left(\frac{1}{h(k)}+k^T A(A+B_k)^{-1}B_kk\right).
\end{multline}

Finally, by applying the Theorem~4.1.11 in Dembo and Zeitouni~\cite{dembo2009large}, we conclude the strong moderate deviation principle.

\end{proof}

\subsection{Proofs for the total number of connected components.}
Now we give the proof of Theorem~\ref{total}.

\begin{proof}
The proof begins by utilizing the symmetry property of the independent and identically distributed sequence $\{X_i^n\}$. For any integer $j \in \{1,\dots,n\}$, any $0<\alpha\le \mu^n$ and any vector $k \in \N^d$ such that $\alpha n <k\le \mu^n n$, the following inequality holds for the conditional probability:
\begin{multline}\label{eq54}
\P\left(N(Z_nn)=j-1\right)\P\left(\sum_{i=1}^{j-1}X_i^n=\mu^n n-k, X_i^n\leq \alpha n, \forall i\right)\frac{Z_nn\P(X^n=k)}{\P\left(\sumn X_i^n=\mu^n n\right)}\\
\leq
\P\left(N(Z_n n)=j,\max_{1\leq i\leq N(nZ_n)}X_i^n=k\bigg| \sumn X_i^n=\mu^nn\right)\\
    \leq \P(N(Z_n n)=j-1)\P\left(\sum_{i=1}^{j-1}X_i^n=\mu^n n-k\right)\frac{Z_nn\P(X^n=k)}{\P\left(\sumn X_i^n=\mu^n n\right)}.
\end{multline}
\paragraph{\bf Exponential Tightness.}
Thanks to Proposition \ref{asymp} and the moderate deviation principle of Poisson random variables, for any $\xi>0$, we obtain the following upper bound
\begin{multline*}
    \mdpsup \left(
    \frac{1}{a_n\sqrt{n}}\bigg|N(Z_nn)-\left(|\cc|-\frac{1}{2}\langle\cc,\kappa\cc\rangle\right)n\bigg|>\xi,A_n(M)\bigg|\sumn X_i^n=\mu^n n
    \right)\\
    \leq\mdpsup \left(
    \frac{1}{a_n\sqrt{n}}\bigg|N(Z_nn)+1-
    \left(|\cc|-\frac{1}{2}\langle\cc,\kappa\cc\rangle\right)n\bigg|>\xi\right)\\
  -\frac{1}{2}\inf_{|y|\le M}\langle D^{-1}_{\mu-\cc}(I-D_{\cc}\kappa)y,(I-D_{\mu}\kappa)y\rangle\\
    \leq -\frac{\xi^2}{2}\frac{1}{|\cc|-\frac{1}{2}\langle\cc,\kappa\cc\rangle}-\frac{M^2}{2}\lambda_2,
\end{multline*}
where $\lambda_2$ denotes the minimum eigenvalue of the symmetric matrix $(I-\kappa D_{\mu})D^{-1}_{\mu-\cc}(I-D_c\kappa)$. 
By combining this result with Lemma~\ref{et1} and choosing $M$ sufficiently large, , followed by taking the limit $\xi\rightarrow\infty$, we establish the exponential tightness for the sequence of conditional distributions of $N(Z_n n)$.

\paragraph{\bf Upper bound.}To prove the upper bound of the MDP, we apply equation \eqref{mdpsum} and Chebyshev's inequality. This yields that for all $z\in\R^d$, 
\begin{multline*}
    \frac{1}{a_n^2}\log \P\left(\sum_{i=1}^{\lfloor(|\cc|-\frac{1}{2}\langle\cc,\kappa\cc\rangle )n+xa_n\sqrt{n}-1\rfloor}X_i^n=\lceil\cc^n n-ya_n\sqrt{n}\rceil\right)\\
    \leq\frac12\bigg\langle z,\left(\Phi-\frac{\cc\cc^T}{|\cc|-\frac12\langle\cc,\kappa\cc\rangle}\right)z\bigg\rangle
    +\langle z,y\rangle+\frac{2x}{2|\cc|-\langle\cc,\kappa\cc\rangle}\langle z,\cc\rangle
    +O_z\left(\frac{1}{a_n}+\frac{a_n}{\sqrt{n}}\right).
\end{multline*}
By choosing the optimal parameter $z$ (through a standard Legendre-Fenchel transformation argument),
\[
z_y=-\left(\Phi-\frac{\cc \cc^T}{|\cc|-\frac 12\langle \cc, \kappa\cc\rangle}\right)^{-1}\left(\frac{\cc x}{|\cc|-\frac 12\langle \cc, \kappa\cc\rangle}+y\right),
\]
and using Proposition \ref{asymp} and upper bound in ~\eqref{eq54}, we derive the following upper bound for the joint conditional probability
\begin{multline}
    \label{uptotal}
    \sup_{x\in\overline\B}\sup_{|y|\le M}\frac{1}{a_n^2}\log\P\bigg(N(Z_nn)=\left(|\cc|-\frac{1}{2}\langle \cc,\kappa\cc\rangle \right)n+a_n\sqrt{n}x,\\ 
    \max_{1\leq i\leq N(nZ_n)}X_i^n=(\mu^n-c^n)n+a_n\sqrt{n}y\bigg|\sumn X_i^n=\mu^n n\bigg)\\
    \leq \sup_{x\in\overline\B}\mdp\left(N(Z_nn)=\left(|\cc|-\frac{1}{2}\langle \cc,\kappa\cc\rangle \right)n+a_n\sqrt{n}x-1\right)+\\  
    \sup_{x\in\overline\B}\sup_{|y|\le M}\bigg(-
    \frac{1}{2}\langle
    D_{\mu-\cc}^{-1}(I-D_{\cc} \kappa)y,(I-{D_{\mu}} \kappa)y
    \rangle\\    
    + \frac12
    \langle z_y, \left(\Phi-\frac{\cc \cc^T}{|\cc|-\frac 12\langle \cc, \kappa\cc\rangle}\right) z_y\rangle
    +\langle z_y,y\rangle+\frac{2x}{2|\cc|-\langle\cc,\kappa\cc\rangle}\langle z_y,\cc\rangle\bigg)    \\
    +O_{M,\beta}\left(\frac{1}{a_n}+\frac{a_n}{\sqrt{n}}\right)+O_M\left(\frac{1}{a_n^2}+\frac{a_n}{\sqrt{n}}\right).
\end{multline}
Recalling the definitions of the matrices $A$ and $B$ from Theorem~\ref{thm2} and Theorem~\ref{thm3},  the upper bound above simplifies as follows:
\begin{multline}\label{ss}
\sup_{x\in\overline\B}
    \sup_{|y|\le M}
    \bigg(-\frac{x^2}{2}
    \frac{1}{|\cc|-\frac{1}{2}\langle\cc,\kappa\cc\rangle}\\
    -\frac{1}{2}y^TAy
-\frac12 \left(\frac{\cc x}{|\cc|-\frac 12\langle \cc, \kappa\cc\rangle}+y\right)^TB\left(\frac{\cc x}{|\cc|-\frac 12\langle \cc, \kappa\cc\rangle}+y\right)\bigg).
\end{multline}
Thanks to Lemma~\ref{pdefinite2}, the matrix 
$A+B$
is positive definite.  For $M$ sufficiently large,  the internal extremum over $y$ is obtained at 
\[
y_x=-(A+B)^{-1}B\frac{\cc x}{|\cc|-\frac 12\langle \cc, \kappa\cc\rangle}.
\]
Substituting $y_x$ into the upper bound~\eqref{ss},  the supremum becomes
\[
\sup_{x\in\overline\B}\left\{-\frac{x^2}{2}
    \frac{1}{|\cc|-\frac{1}{2}\langle\cc,\kappa\cc\rangle}-\frac12\left(\frac{x}{|\cc|-\frac 12\langle \cc, \kappa\cc\rangle}\right)^2
    \cc^TA(A+B)^{-1}B\cc\right\}.
\]
Using the fact
\[
\text{Card}\left\{x\in\B\bigg|(|\cc|-\frac{1}{2}\langle\cc,\kappa\cc\rangle )n+xa_n\sqrt{n}-1\in\{1,\cdots n\}\right\}\leq n,
\]
and the upper bound in the relation~\eqref{uptotal} with Lemma~\ref{et1}, for $M$ sufficiently large, we establish the  upper bound for the MDP:
\begin{multline}
\lim_{\delta\rightarrow 0}\mdpsup\left(\frac{1}{a_n\sqrt{n}}
\left(N(Z_n n)-\left(|\cc|-\frac{1}{2}\langle\cc,\kappa\cc\rangle\right)n\right)\in\overline{\B}\bigg|\sumn X_i^n=\mu^n n\right)\\
\leq -\frac{\beta^2}{2}
    \frac{1}{|\cc|-\frac{1}{2}\langle\cc,\kappa\cc\rangle}-\frac12\left(\frac{\beta}{|\cc|-\frac 12\langle \cc, \kappa\cc\rangle}\right)^2
    \cc^TA(A+B)^{-1}B\cc. 
\end{multline}
\paragraph{\bf Lower bound}
To prove the lower bound, we first fix $y_\beta$ as determined in the proof of upper bound.
Now for  any $\eps>0$ and $n$ sufficiently large, we have the following lower bound,
\begin{multline*}
  \P\left(\frac{1}{a_n\sqrt{n}}\left(N(Z_n n)-\left(|\cc|-\frac{1}{2}\langle\cc,\kappa\cc\rangle\right)n\right)\in\B\bigg|\sumn X_i^n=\mu^n n\right)\\
  \ge
    \P\Bigg(N(Z_n n)=\bigg\lfloor\left(|\cc|-\frac{1}{2}\langle\cc,\kappa\cc\rangle\right)n+a_n\sqrt{n}\beta\bigg\rfloor,\\\max_{1\leq i\leq N(nZ_n)}X_i^n=k^n(\eta),\eta\in B_\eps(y_\beta)\bigg|\sumn X_i^n=\mu^n n\Bigg).
\end{multline*}
Hence, using the lower bound in~\eqref{eq54} Proposition \ref{asymp}, Proposition \ref{MDPcppp}, Lemma \ref{term3} and Stirling's formula, we obtain the lower bound
\[
\begin{aligned}
    &\lim_{\delta\rightarrow 0}\mdpinf \left(\frac{1}{a_n\sqrt{n}}\left(N(Z_n n)-\left(|\cc|-\frac{1}{2}\langle\cc,\kappa\cc\rangle\right)n\right)\in\B\bigg|\sumn X_i^n=\mu^n n\right)\\
    \geq&\mdpinf\left(
    N(Z_n n)=\bigg\lfloor\left(|\cc|-\frac{1}{2}\langle\cc,\kappa\cc\rangle\right)n+a_n\sqrt{n}\beta\bigg\rfloor-1
    \right)\\
    &\qquad\qquad\qquad+\lim_{\eps\rightarrow 0}
    \liminf_{n\rightarrow\infty}\inf_{\eta\in B_{\eps}(y_\beta)}\frac{1}{a_n^2}\log \frac{Z_nn\P\left(X^n=k^n(\eta)\right)}{\P\left(\sumn X_i^n=\mu^n n\right)}\\
    +&\lim_{\eps\rightarrow 0}\mdpinf\left(
    \frac{1}{a_n \sqrt{n}}\left(\sum_{i=1}^{\lfloor\left(|\cc|-\frac{1}{2}\langle\cc,\kappa\cc\rangle\right)n+a_n\sqrt{n}\beta\rfloor-1}X_i^n-\cc^n n\right)\in B_{\eps}(-y_\beta)\bigg| X_i^n\leq \alpha n
    \right)\\
    \geq&
    -\frac{\beta^2}{2}\frac{1}{|\cc|-\frac{1}{2}\langle\cc,\kappa\cc\rangle}
    - \frac{1}{2}\langle
    D_{\mu-\cc}^{-1}(I-D_{\cc}\kappa)y_\beta,(I-D_{\mu}\kappa)y_\beta\rangle\\
    &\qquad\qquad-\frac12 \left(\frac{\cc \beta}{|\cc|-\frac 12\langle \cc, \kappa\cc\rangle}+y_\beta\right)^T\left(\Phi-\frac{\cc \cc^T}{|\cc|-\frac 12\langle \cc, \kappa\cc\rangle}\right)^{-1}\left(\frac{\cc \beta}{|\cc|-\frac 12\langle \cc, \kappa\cc\rangle}+y_\beta\right)\\
    &=-\frac{\beta^2}{2}
    \frac{1}{|\cc|-\frac{1}{2}\langle\cc,\kappa\cc\rangle}-\frac12\left(\frac{\beta}{|\cc|-\frac 12\langle \cc, \kappa\cc\rangle}\right)^2
    \cc^TA(A+B)^{-1}B\cc.
\end{aligned}
\]

Finally, by applying the Theorem~4.1.11 in Dembo and Zeitouni~\cite{dembo2009large}, we establish the strong moderate deviation principle.

\end{proof}

%% file: app.tex
\section{The proofs for the connection probability}\label{appen}

We first show the proof of the upper bound.
\begin{proof}[Proof of Lemma~\ref{esti2p}]
As proved in Lemma 4.10 in~\cite{LDPinhomo}, given a $k\in\N^d$ with $k_r\geq1$ for some $r\in\S$ and for any $m\in\N^d$ with $m\geq k$, we have
    \begin{equation}\label{4.10cite}
        p_n(k)\leq \left[\prod_{s\in\S}\binom{m_s-\delta_{r,s}}{k_s-\delta_{r,s}}\right]^{-1}\prod_{s,s'\in \mathcal{S}}\left(1-\frac{\kappa(s,s')}{n}\right)^{-k_s(m_{s'}-k_{s'})},
    \end{equation}
    where
    \[
        \delta_{r,s}=\left\{
        \begin{aligned}
            &1, r=s,\\
            &0,r\neq s.
        \end{aligned}
        \right  .
    \]
Let $m^n_r=\lfloor k_r^n(1-e^{-\frac{1}{n}(\kappa k^n)_r})^{-1}\rfloor $ be the chosen sequence in \eqref{4.10cite}. Applying the Stirling's formula $n!=\sqrt{2\pi n}(\frac{n}{e})^n(1+O(\frac{1}{n}))$, we have
 \begin{multline*}
\left[\prod_{s\in\mathcal{S}}\binom{m^n_s-\delta_{r,s}}{k^n_s-\delta_{r,s}}\right]\prod_{s,s'\in \mathcal{S}}\left(1-\frac{\kappa(s,s')}{n}\right)^{k^n_s(m^n_{s'}-k^n_{s'})}
\\
=\prod_{s\in\S}\left(\frac{1}{\sqrt{2\pi k^n_s}}\right){\exp\left(-\frac{1}{2n}\sum_{s\in\mathcal{S}}(\kappa k^n)_s\right)}
\bigg/
{\prod_{s\in\mathcal{S}}\left(1-e^{-\frac{1}{n}(\kappa k^n)_s}\right)^{k^n_s}}\left(1+O\left(\frac{1}{n}\right)\right).
\end{multline*}
The upper bound is then obtained directly by applying Lemma~\ref{lm3} to this asymptotic expression.
\end{proof}

Now we give the proof of the lower bound.
\begin{proof}[Proof of Lemma~\ref{p_lb}]
    Let $y=\mu-c$. 
    For any $\delta$ such that $0< \delta<\frac13 \inf\{y_s:s\in\S\}$, we define $y^{(\delta)}:=(y_s-\delta,s\in\cal{S})$ and 
    \[
    \nu^{(\delta)}:=\left(\frac{y^{(\delta)}_s}{1-e^{-(\kappa y^{(\delta)})_s}},  s\in\S\right).
    \]
    It is clear that $\nu^{(\delta)}$ converges to $\nu^{(0)}=:\nu$ as $\delta\rightarrow0$. For all $\beta\in\R$, let $m^{(n,\delta)}(\beta)\in\N^d\backslash\{0\}$ be such that $m_s^{(n,\delta)}(\beta)=\lfloor(k_s^n(\beta)-n\delta)(1-e^{-\frac 1n (\kappa(k^n(\beta)-n\delta))_s})^{-1}\rfloor$ for all $s\in\cal{S}$. Now fixing $r\in\S$ with $y_r>0$, for $l\in\N^d\backslash\{0\}$ satisfies $l\leq m^{(n,\delta)}$ and $l_r\ge 1$, we define
    \[
    b_n^{(\delta)}(l)(\beta):=p_n(l)\left[\prod_{s\in\mathcal{S}}\binom{m^{(n,\delta)}_s(\beta)-\delta_{r,s}}{l_s-\delta_{r,s}}\right]\prod_{s,s'\in \mathcal{S}}\left(1-\frac{\kappa(s,s')}{n}\right)^{l_s(m^{(n,\delta)}_{s'}(\beta)-l_{s'})}.
    \]
     Using Stirling's approximation, for any $l\leq m^{(n,\delta)}(\beta)$, $l_r\geq 1$ and $|l|>\epsilon n$ with some $\epsilon>0$, we have the following.
    \begin{equation}\label{stirling}
    \begin{aligned}
    b_n^{(\delta)}(l)(\beta)
    &\leq\prod_{s\in\S}\sqrt{\frac{m_s^{(n,\delta)}(\beta)-\delta_{r,s}}{2\pi (l_s-\delta_{r,s})(m_s^{(n,\delta)}(\beta)-l_s)}}\left(1+O\left(\frac1n\right)\right)e^{-nf(\frac{l}{n},\frac{m^{(n,\delta)}(\beta)}{n})}p_n(l),
    \end{aligned}
    \end{equation}
    where $f$ is defined in equation (4.34) of \cite{LDPinhomo}
    \[
    f(x,x'):=\<x,\log\frac{x}{x'}\>+\<x'-x,\log\frac{x'-x}{x'}\>+\<x'-x, \kappa x\>.
    \] 
  In order to obtain an upper bound of~\eqref{stirling}, we first notice that for $n$ sufficiently large and $\delta$ sufficiently small, we have
    \begin{equation*}
    \frac{m_s^{(n,\delta)}(\beta)}{m_s^{(n,\delta)}(\beta)-k^n_s(\beta)}
    \leq \frac{3|\mu-c|}{e^{-\Vert\kappa\Vert_\infty\vert\mu-c\vert}\min_{r\in \S}(\mu_r-c_r)}.
    \end{equation*}
    Hence, using a similar approach to the proof of Claim 3 in Lemma 4.12 in \cite{LDPinhomo}, 
    we know that for any $l\in[k^n(\beta)-2n\delta, k^n(\beta)]$, there exists a positive constant $C_0$, such that
    \begin{equation}\label{l in}
    \begin{aligned}
    b_n^{(\delta)}(l)(\beta)
    \leq\frac{C_0
     }{\prod_{s\in\S}\sqrt{k^n_s(\beta)-2n\delta}}\left(1+O\left(\frac1n\right)\right)e^{nC^n(\delta)}e^{-nf\left(\frac{k^n(\beta)}{n},\frac{m^{(n)}(\beta)}{n}\right)}p_n(k^n(\beta)),
    \end{aligned}
    \end{equation}
     where
     \begin{multline*}
     C^n(\delta)=
     \sup_{|\beta|\leq M}\sup_{l\in[k^n(\beta)-2n\delta,k^n(\beta)]}\Bigg|f\left(\frac{k^n(\beta)}{n},\frac{m^{(n)}(\beta)}{n}\right)-
     \\f\left(\frac{l}{n},\frac{m^{(n,\delta)}(\beta)}{n}\right)-\left\langle \frac{k^n(\beta)-l}{n},\log\left(1-e^{-\frac{1}{n}\kappa l}\right)\right\rangle\Bigg|.
     \end{multline*}
     It is easy to see that $\Bigg|\left(\frac{k^n(\beta)}{n},\frac{m^{(n)}(\beta)}{n}\right)-\left(\frac{l}{n},\frac{m^{(n,\delta)}(\beta)}{n}\right)\Bigg|=O(\delta)$ uniformly for $|\beta|\leq M$,  $n$ sufficiently large, and $l\in[k^n(\beta)-2n\delta,k^n(\beta)]$. Since $f(\cdot,\cdot)$ is smooth at $(y,\nu)$, we have $C^n(\delta)= O(\delta)$ uniformly for $n$ sufficiently large. 
     
     When $l\notin[k^n(\beta)-2n\delta, k^n(\beta)]$ and $|l|>\epsilon n$, we use Stirling's approximation \eqref{stirling} with the results from Lemma \ref{esti2p} to obtain that, for $n$ sufficiently large,
     \[
     \sup_{|\beta|\leq M}b_n^{(\delta)}(l)(\beta)\leq O(n^d)e^{-n\tilde{C}_n(\delta)}, 
     \]
     where 
     \[\begin{aligned}
     \tilde{C}_n(\delta)=\inf_{x\in\cal{U}_n} \sum_{s\in\S}\left(x_s\log\frac{\frac{nx_s}{m^{(n,\delta)}_s}}{1-e^{-(\kappa x)_s}}+(\nu^{(\delta)}_s-x_s)\log\frac{1-\frac{nx_s}{m^{(n,\delta)}_s}}{e^{-(\kappa x)_s}}\right),
     \end{aligned}
     \]
     with
     \[
\cal{U}_n:=\left\{x\in\R_+^d:{|x|>\epsilon, \left|x-\left(\frac{k^n}{n}-\delta\right)\right|>\delta},x\leq\frac{m^{(n,\delta)}}{n}\right\}.
     \]
     We can rewrite the right-hand side of $\tilde{C}_n(\delta)$ as
\[
\inf_{x\in\cal{U}_n}\sum_{s\in\S}\<\nu_s^{(\delta)},H\left(Q^{(n,\delta)}_s(x),T_s(x)\right)\>,
\]
where  $H(q,t)=q\log\frac{q}{t}+(1-q)\log\frac{1-q}{1-t}$, $Q^{(n,\delta)}_s(x)=nx_s\backslash m^{(n,\delta)}_s$ and $T_s(x)=1-e^{-(\kappa x)_s}$.
We can observe that  the function
     $H(q,t)$ 
     is structured as a Kullback-Leibler divergence between two Bernoulli distributions with parameters $q$ and $t$. Therefore, $H(q,t)$ is non-negative for all $q,t\in[0,1]$ and the equality $H(Q^{(n,\delta)}_s(x),T_s(x))=0$ is maintained if and only if $Q^{(n,\delta)}_s(x)=T_s(x)$ for all $s\in\S$. This condition is satisfied precisely when $x=\frac{k^n}{n}-\delta$ or $x=0$. Let $J^{(n,\delta)}$  be the Jacobian matrix of $T(x)-Q^{(n,\delta)}(x)$ for $x=\frac{k^n}{n}-\delta$ and $J=D^{-1}_{\mu}(I-D_c\kappa )$. 
     By a direct computation combined with the dual equation, we get
     \[
     J^{(n,\delta)}_{sr}=
     \begin{cases}
         \frac{n}{m_s^{(n,\delta)}}-\kappa_{ss}e^{-\left(\kappa(\frac{k^n}{n}-\delta)\right)_s}&s=r,\\
         \kappa_{rs}e^{-\left(\kappa(\frac{k^n}{n}-\delta)\right)_s}&s\neq r.
     \end{cases}
     \]
     and $     J=\lim_{\substack{n\to \infty\\
     \delta\to 0}}J^{(n,\delta)}$.
Since the smallest eigenvalue of $J$ is greater than zero, then for $\delta$ sufficiently small and $n$ sufficiently large, there exists $C_1>0$ such that 
     \[
     \inf_{x\in\cal{U}_n}\frac{x^T(J^{(n, \delta)})^2x}{\|x\|^2}>C_1.
     \]
     By the equivalence of the norms $\ell_1$ and $\ell_2$ in $\R^d$, there exists $C_2>0$ such that $|T(x)-Q^{(n,\delta)}(x)|>C_2\delta$ whenever $|x|>\epsilon$ and $|x-\frac{k^n}{n}+\delta|>\delta$. Since $\nu^{(\delta)}_s>0$ for all $s\in\S$ and the Hessian matrix of $H(\cdot, \cdot)$ is positive definite when $x=y^{(\delta)}$, there exists $C_2>0$ such that $\<\nu^{(\delta)},H(x,\nu^{(\delta)})\>>C_2\delta^2$.     
  Using the facts $\#\{l\in\N^d:l\leq m^{(n,\delta)}\}= O(n^d)$ and $f\left(\frac{k^n}{n},\frac{m^{(n)}}{n}\right)=\exp\left(-\langle k^n,\log (1-e^{-\kappa \frac{k^n}{n}})\rangle\right)$,  for $\frac{1}{\sqrt{n}}
     \ll\delta_n\ll\frac{a_n^2}{n}$ and any sufficiently small $\epsilon>0$, we have
     \[
     \eps_n:=\sup_{|\beta|\leq M}\sum_{\substack{ l:|l|>\epsilon n,\\ l\notin[k^n(\beta)-2n\delta_n, k^n(\beta)]}} b_n^{(\delta_n)}(l)(\beta)\rightarrow 0 \text{ as } n\rightarrow\infty.
     \]
Moreover, for all $\beta$ such that $|\beta|\leq M$, thanks to  the inequality \eqref{l in}, we obtain
    \begin{multline*}
    \sum_{l\in[k^n(\beta)-2n\delta_n, k^n(\beta)]}b_n^{(\delta_n)}(l)(\beta)\leq C_\kappa \prod_{s\in\S}\left(\frac{n}{\sqrt{k^n_s(\beta)-2n\delta_n}}\right)\left(1+O\left(\frac1n\right)\right)e^{nC^n(\delta_n)}\\
    \exp\left(-\langle k^n(\beta),\log (1-e^{-\kappa \frac{k^n(\beta)}{n}})\rangle\right)p_n(k^n(\beta)),
    \end{multline*}
    where $C_{\kappa}$ is a positive constant that depends only on $\kappa$. Thus, for all $\beta$ such that $|\beta|\leq M$, the sum of $b_n^{(\delta_n)}(l)(\beta)$ for all $l$ with $|l|>\epsilon n$ has the following upper bound:
 \begin{equation}\label{citeclaim3}
    \begin{aligned}
    \sum_{l:|l|>\epsilon n}b_n^{(\delta_n)}(l)(\beta)\leq C_\kappa &\prod_{s\in\S}\left(\frac{n}{\sqrt{k^n_s(\beta)}}\right)\left(1+O\left(\frac1n\right)\right)\\
    &\exp\left(-\langle k^n(\beta),\log (1-e^{-\kappa \frac{k^n(\beta)}{n}})\rangle\right)e^{nC^n(\delta)}p_n(k^n(\beta))+\eps_n,
    \end{aligned}
    \end{equation}
    where $\eps_n\rightarrow0$ as $n\rightarrow\infty$ and $C^n(\delta)=o\left(\frac{a_n^2}{n}\right)$. 
    
    For the term $\sum_{l:|l|\leq \lfloor a_n\rfloor}b_n^{(\delta)}(l)(\beta)$, we claim that
  \begin{equation}\label{citeclaim1}
 \limsup_{n\to\infty}\sup_{|\beta|\leq M}\sum_{l:|l|\leq \lfloor a_n\rfloor}b_n^{(\delta_n)}(l)(\beta)\leq 1-\frac{y_r}{\nu_r}.
    \end{equation}
To show this limit, we first apply \eqref{estimate} in Lemma \ref{est} and get
   \begin{equation}\label{claim21}
   \begin{aligned}
       &\sum_{l:|l|\leq \lfloor a_n\rfloor}b_n^{(\delta_n)}(l)(\beta)\\
       \leq&\frac{n}{m^{(n,\delta_n)}_r(\beta)}\sum_{l:|l|\leq \lfloor a_n\rfloor}\tau(l)l_r\prod_{s\in\S}\left[\frac{1}{l_s!}\left(\frac{m^{(n,\delta_n)}_s(\beta)}{n}e^{-\left(\kappa\frac{ m^{(n,\delta_n)}(\beta)-l}{n}\right)_s}\right)^{l_s}\right].
   \end{aligned}
   \end{equation}
To estimate the upper bound, we need the following two relations,
 \begin{equation}\label{claim23}
    \begin{aligned}
    &\sum_{l:|l|\leq \lfloor a_n\rfloor}\tau(l)l_r\frac{1}{\prod_{s\in\S}l_s!}\left|\prod_{s\in\S}\left(\frac{m^{(n,\delta_n)}(\beta)}{n}e^{-(\kappa\frac{ m^{(n,\delta_n)}(\beta)-l}{n})_s}\right)^{l_s}-\prod_{s\in\S}\left(\nu_se^{-(\kappa \nu)_s}\right)^{l_s}\right|\\
    &= \sum_{l:|l|\leq \lfloor a_n\rfloor}\tau(l)l_r\frac{1}{\prod_{s\in\S}l_s!}\prod_{s\in\S}\left(\nu_se^{-(\kappa \nu)_s}\right)^{l_s}\\
    &\qquad\qquad\qquad\qquad\left|\prod_{s\in\S}\left(\frac{m_s^{(n,\delta_n)}(\beta)}{n\nu_s}e^{-\left(\kappa\frac{ m^{(n,\delta_n)}(\beta)-l}{n}\right)_s+(\kappa\nu)_s}\right)^{l_s}-1\right|,
    \end{aligned}
    \end{equation}
and
    \begin{equation}\label{claim22}
    \begin{aligned}
    &\left|\log\prod_{s\in\S}\left(\frac{m_s^{(n,\delta_n)}(\beta)}{n\nu_s}e^{-\left(\kappa\frac{ m^{(n,\delta_n)}(\beta)-l}{n}\right)_s+(\kappa\nu)_s}\right)^{l_s}\right|\\
    =&\left|\sum_{s\in\S}l_s\log\frac{m_s^{(n,\delta_n)}(\beta)}{n\nu_s}-l^T\kappa\left(\frac{m^{(n,\delta_n)}(\beta)}{n}-\nu\right)+\frac{l^T\kappa l}{n}\right|.
    \end{aligned}
    \end{equation}
    For $s\in\S$ and $x\in\R^d$, define $g_s(x)=x_s(1-e^{-(\kappa x)_s})^{-1}$. Since $\left|\frac{k^n(\beta)}{n}-y\right|=O_\beta\left(\frac{a_n}{\sqrt{n}}\right)$ and the function $g$ is continuous,  the terms in the right-hand side of equation \eqref{claim22} have the following upper bounds,
\[\begin{aligned}
\Bigg|\kappa&\left(\frac{m^{(n,\delta_n)}(\beta)}{n}-\nu\right)\Bigg|
\le 
d\|\kappa\|_\infty\left(\left|\frac{m^{(n,\delta_n)}(\beta)}{n}-\nu^{(\delta_n)}\right|
+
|\nu^{(\delta_n)}-\nu|\right)\\
&=d\|\kappa\|_\infty\left(\left|g\left(\frac{k^n(\beta)}{n}-\delta_n\right)- g(y-\delta_n)\right|
+
|g\left(y-\delta_n\right)- g(y)|\right)
,
\end{aligned}\]
and
\[
\begin{aligned}
\Bigg|\log&\left(\frac{m_s^{(n,\delta_n)}(\beta)}{n\nu_s}\right)\Bigg|\leq \left|\log\frac{m_s^{(n,\delta_n)}(\beta)}{n}-\log\nu_s^{(\delta_n)}\right|
+
|\log\nu_s^{(\delta_n)}-\log\nu_s|\\
&=\left|\log g_s\left(\frac{k^n(\beta)}{n}-\delta_n\right)- \log g_s(y-\delta_n)\right|
+
|\log g_s\left(y-\delta_n\right)- \log g_s(y)|.
\end{aligned}\]
   These two quantities   converge to $0$ as $n\to \infty$ uniformly for $\beta$ such that $|\beta|\leq M$. In conclusion, the difference~\eqref{claim23} goes to $0$ as $n\to\infty$ uniformly for $|\beta|\le M$. Recall that, in Lemma \ref{cor}, we have proved that
    \[
    \sum_{l\in\N^d}\tau(l)l_r\frac{1}{\prod_{s\in\S}l_s!}\prod_{s\in\S}\left(\nu_se^{-(\kappa \nu)_s}\right)^{l_s}=\frac{\nu_r-y_r}{\nu_r},
    \]        
 Thanks to~\eqref{claim21}, we obtain claim \eqref{citeclaim1}.

   Recall that in Claim 2 of Lemma 4.12 in \cite{LDPinhomo}, for some sufficiently small  $\delta_1\in(0,1)$ and sufficiently large $N_1\in\N$, we have 
   \begin{equation}\label{citeclaim2}
    \limsup_{R\rightarrow\infty} \sup_{N>N_1}\sup_{\delta\in(0,\delta_1]}\sup_{|\beta|\leq M}\sum_{R<|l|\leq\epsilon N}b_n^{(\delta)}(l)(\beta)=0,
    \end{equation}
which implies
\begin{equation}\label{citeclaim2}
    \limsup_{n\to \infty}\sup_{|\beta|\leq M}\sum_{\lfloor a_n\rfloor<|l|\leq\epsilon n}b_n^{(\delta_n)}(l)(\beta)=0.
    \end{equation}
By combining \eqref{citeclaim1} and \eqref{citeclaim2} and taking $\epsilon^*=\frac{y_r}{4\nu_r}$, for $n$ sufficiently large, we have
\begin{equation}\label{<epsion n}
    \sup_{|\beta|\leq M}\sum_{l:|l|\leq\epsilon n}b_n^{(\delta_n)}(l)(\beta)\leq 1-\frac{y_r}{\nu_r}+\epsilon^*=1-\frac{3y_r}{4\nu_r}.
\end{equation}
We also know from Lemma 4.9 in \cite{LDPinhomo}, 
    \begin{equation}\label{lemma4.9}
    \sum_{l\in \N^d\backslash\{0\}:l \leq m^{(n,\delta_n),l_r\ge 1}}b_n^{(\delta_n)}(l)(\beta)=1.
    \end{equation}
    We combine \eqref{<epsion n} and \eqref{lemma4.9} to deduce that for $n$ sufficiently large,
    \[
    \inf_{|\beta|\leq M}\sum_{|l|> \eps n}b_n^{(\delta_n)}(l)(\beta)
    >\frac{3y_r}{4\nu_r}>0,
    \]
    which is bounded away from $0$. Using the inequality \eqref{citeclaim3}, we have
    \begin{multline*}
    \inf_{|\beta|\leq M}C_\kappa \prod_{s\in\S}\left(\frac{n}{\sqrt{k^n_s(\beta)}}\right)\left(1+O\left(\frac1n\right)\right)\\
    \exp\left(-\langle k^n(\beta),\log (1-e^{-\kappa \frac{k^n(\beta)}{n}})\rangle\right)e^{-nC^n(\delta_n)}p_n(k^n(\beta))+\eps_n\\
   \geq
    \frac{3y_r}{4\nu_r}>0.
    \end{multline*}
  Since $\eps_n\rightarrow0 $ as $n\rightarrow\infty$, we can take $n$ sufficiently large such that  $\eps_n<\frac{y_r}{4\nu_r}$, and get
 \begin{multline*}
      \inf_{|\beta|<M}C_\kappa \prod_{s\in\S}\left(\frac{n}{\sqrt{k^n_s(\beta)}}\right)\left(1+O\left(\frac1n\right)\right)\exp\left(-\langle k^n(\beta),\log (1-e^{-\kappa \frac{k^n(\beta)}{n}})\rangle\right)e^{nC^n(\delta_n)}p_n(k^n(\beta))\\
  \geq
    \frac{y_r}{2\nu_r}>0.
\end{multline*}
Taking logarithms on both sides of the above inequality, dividing by $a_n^2$ and taking the limit, we obtain the following
\[
        \liminf_{n\to\infty}\inf_{\beta:|\beta|\leq M}\frac{1}{a_n^2}\log \prod_{s\in\S}\left(\frac{n}{\sqrt{k^n_s(\beta)}}\right)\exp\left(-\langle k^n(\beta),\log (1-e^{-\kappa \frac{k^n(\beta)}{n}})\rangle\right)p_n(k^n(\beta))\geq 0.
\]
\end{proof}